%% file: M144477_R2.tex
\documentclass[final,onefignum,onetabnum]{siamart171218}

\usepackage{amsmath,latexsym,amssymb}

\usepackage[cp1251]{inputenc} 

  \usepackage{multicol}
\usepackage{fancybox}
     \usepackage{xcolor}
\usepackage{amsmath}
\usepackage{cite}

  \usepackage{enumerate}

\newtheorem{example}{Example}
\renewcommand{\baselinestretch}{1.0}
\def\eqsp{\noalign{\vskip 5pt}}
\def\dsp{\displaystyle}

\usepackage{bm}
\usepackage{graphicx}
\usepackage{framed}

\usepackage[all]{xy}

\input{zmacros.tex}

\usepackage{multicol}
\usepackage{fancybox}
     \usepackage{xcolor}

\usepackage{cite}

\usepackage{cases}
\usepackage{multirow}
\usepackage{multicol}
\usepackage{arydshln}
\usepackage{amsmath}
\usepackage{amssymb}
\usepackage{epstopdf}
\usepackage{enumitem}

\renewcommand{\baselinestretch}{1.0}

\usepackage{bm}
\usepackage{framed}

\usepackage[all]{xy}
\usepackage{subfigure}

\newcommand{\bfalf}{\mbox{\boldmath $\alpha$}}
\newcommand{\bfbeta}{\mbox{\boldmath $\beta$}}
\newcommand{\bfgamma}{\mbox{\boldmath $\gamma$}}

\graphicspath{{Figs/}}

\begin{document}

\title{High order compact schemes for  flux  type~BCs \thanks{Submitted to the journal's Methods and Algorithms for Scientific Computing section.
\funding{Z. Li is partially supported by a Simon's grant 633724. K. Pan is supported by the National Natural Science Foundation of China (No. 41874086) and the Excellent Youth Foundation of Hunan Province of China (No. 2018JJ1042).}}}


\author{Zhilin Li \thanks{CRSC and Department of Mathematics, North Carolina State University, Raleigh, NC 27695-8205, USA (\email{zhilin@math.ncsu.edu}).}
\and
Kejia Pan\thanks{School of Mathematics and Statistics, HNP-LAMA, Central South University, Changsha  410083, Hunan, China (\email{kejiapan@csu.edu.cn}).}
}

\date{}
\maketitle


\begin{abstract}
In this paper new innovative fourth order compact schemes for Robin and Neumann boundary conditions have been developed  for boundary value problems of elliptic PDEs in two and three dimensions.
 Different from traditional finite difference operator approach, which may not work for flux type of boundary conditions,  carefully designed undetermined coefficient methods  are utilized in developing high order compact   (HOC) schemes.
The new methods  not only can be utilized to design HOC schemes for flux type of boundary conditions but can also be applied to general elliptic PDEs including Poisson, Helmholtz, diffusion-advection, and anisotropic equations with linear boundary conditions.
 In the new developed HOC methods, the coefficient matrices are generally M-matrices,  which guarantee the discrete maximum principle for well-posed problems, so  the convergence of the HOC methods. The developed HOC methods are  versatile and can cover most of high order compact schemes in the literature. The HOC methods for Robin boundary conditions and for anisotropic diffusion and advection equations  with Robin or even
Dirichlet boundary conditions are  likely the first ones that have ever been developed. With the help of pseudo-inverse, or SVD solutions, we have also observed that the developed HOC methods usually have  smaller error constants compared with traditional HOC methods when applicable.  Non-trivial  examples with large wave numbers and oscillatory solutions are presented to confirm the performance of the new HOC methods.
 \end{abstract}



\section{Introduction}
The original motivation of this research is to develop  fourth order compact finite difference schemes for flux boundary conditions (BCs) such as Neumann and Robin for Poisson/Helmholtz equations for which one can apply the standard
fourth order compact scheme at  interior grid points.
After initial success in developing HOC schemes for Robin/Neumann boundary conditions, we have found out that the methodology can be applied to general elliptic boundary value problems (BVP) of the following
\eqm
  &  \grad \cdot \left ( A \grad u\right )  + {\bf a } \cdot \grad u + K u = f(\mathbf{x}), \qquad \mathbf{x}  \in \cal{R},\\
  & \left . \frac{\null}{\null} u \right |_{\partial  { \cal{R}}_1} = u_1(\mathbf{x}), \qquad   \left .  \lp A \frac{\null}{\null} \grad u \cdot \dsp {\bf  n}  + \sigma  u (\mathbf{x})  \rp \frac{\null}{\null} \!\! \right |_{\partial  {\cal{R}}_2} = g(\mathbf{x}), \label{two-BC}
\enm
with a few line changes in the computer codes, where  $\mathbf{x}$ is a point in a rectangular domain $ \cal{R}$,
${\bf n}$ is the unit normal at the boundary pointing outside of the domain. We assume that all the coefficients $A$, ${\bf a }$, $K$, and $\sigma$ are constants, although the methodology developed in this paper can be applied to variable coefficients as well at a cost. We assume that the source term $f(\mathbf{x})\in C^{\nu}$ and the solution $u(\mathbf{x})\in C^{4 +\nu}$ for a $\nu >~0$.  For a fourth order method, we need $\nu=2$ for the convergence proof. The condition $\sigma \ge 0$  is needed to guarantee the well-posedness of the boundary value problem. When $\sigma = 0$, we have a Neumann boundary condition.

One of  advantages of a compact higher order method is that fewer grid points  can be used for the same
order accuracy as a lower order method; therefore, a smaller resulting system of
algebraic equations needs to be solved. This is significant for three or higher dimensional problems to relieve the so-called memory bottleneck. Also, a compact higher order method
maybe needed for highly oscillatory solutions, particularly for wave scattering characterized by Helmholtz equations with large wave numbers $\sqrt{K}$.  Furthermore, high order compact methods are  important for problems with an infinite domain when a mesh size is relatively large.
Moreover,  HOC methods have less grid orientation effects compared with  standard
finite difference schemes since more neighboring grid points in different directions rather than just coordinates ones are involved. On the other hand, it is challenging to develop HOC schemes for flux BCs since there is few information outside of the domain.

One of the earlier  fourth  order compact finite difference (FD) schemes for Poisson equations can be found in \cite{MR612589}. The authors state that the fourth order compact method developed in \cite{MR612589} becomes third order  for Neumann boundary conditions, which can be actually validated in this paper, see Section~\ref{sec:super-3rd}. Other early work can be found in \cite{MR1748146,MR1377177,MR1440337,MR632899} and a few others.  A commonly used  fourth  order compact finite difference scheme for Poisson equations, or Helmholtz equations in 2D  can be found in \cite{strikwerda,li:book-FDFEM,MR1377177}, for example.  For a Helmholtz equation with a Dirichlet boundary condition, one can simply treat the $K u $ term as an external force in the discretization, which does not work for Neumann or Robin boundary conditions, and destroy the M-matrix property for a non-zero $K$.

In deriving HOC schemes,  often  finite difference operators are employed dimension by dimension in discretizing a PDE.   At the boundary of a rectangular domain, Dirichlet boundary have no effect on the accuracy of high order compact schemes. However, for many applications, flux type boundary conditions are provided. For example, in a shear flow, often no-slip (${\bf u}=0$) boundary conditions are prescribed at the top and bottom walls, while at the inlet,  a flux  condition is prescribed. A second order accurate scheme to deal with  flux boundary conditions  would  ruin global fourth order  accuracy if a fourth order discretization is utilized in the interior. Research  on HOC schemes for   flux type boundary conditions can be found,  for example,  in \cite{MR2605462,MR3361665,MR3885818} for time dependent problems based on operator splitting approaches,  for Poisson and
Helmholtz or wave equations \cite{MR4056918,MR4036946,MR3064170,MR2378588,MR2263557,MR2285863},  other related HOC methods and applications \cite{MR2289457,MR2988136,MR3565908,MR2537852,MR4252893,MR4089070,MR3885818,TONG2021107413,MR4156176},
and for diffusion and advection equations \cite{MR1440337,2012Aliao,MR2285863}.
However, few can be found in the literature for general  elliptic partial differential equations with flux boundary conditions. Thus, it is a {\em important problem} to develop HOC  schemes, particularly fourth order, for flux type of boundary conditions (Neumann or  Robin).

There are two important considerations in deriving HOC schemes. The first one is the consistency,  which is relatively easy to confirm. The second one is the stability,  which is more challenging to address and  often left out  in discussions especially for elliptic PDEs in some of research in the literature.  One of tools to ensure the stability is to check whether the coefficient matrix of the system of equations is an M-matrix or not.  The consistence condition plus an M-matrix (stronger stability) condition will  lead to the convergence of the method, see for example, \cite{morton-mayers}.

In this paper, we propose a completely new approach in constructing HOC finite difference schemes for Poisson,  Helmholtz, and diffusion-advection equations with constant coefficients on rectangular domains with Dirichlet, Neumann, or Robin boundary condition  along a part of  boundary of the domain.  The idea is to use linear combinations of the solution  values ($U_{ij}$ or $U_{ijk}$), the source term values ($f_{ij}$ or $f_{ijk}$), and the flux type of boundary condition restricted to the grid points in the compact finite difference stencil. The undetermined coefficients are chosen such that the local truncation errors can be as small as possible in magnitude, $O(h^4)$ in interior  and $O(h^3)$ at the boundary,
while maintaining the discrete maximum principle, that is, the coefficient matrix is an M-matrix when $K\le 0$. We summarize our new HOC methods in this paper below:
\begin{enumerate}[label=\arabic*), itemindent=0.5em]
\item 4th order  compact method for Poisson/Helmholtz equations with Neumann or Robin BCs;
\item 4th order compact  method for diffusion and advection equations with constant coefficients with Dirichlet, Neumann,  or Robin BCs;
\item 4th order compact method for anisotropic diffusion and advection equations with constant coefficients with Dirichlet, Neumann,  or Robin BCs;
\item super 3rd order  compact methods for above equations and BCs without $f$-extensions.
\end{enumerate}
%
In the first three items, we assume that we can extend the source term ($f$) to one grid line (2D) or surface (3D) outside of the boundary with $O(h^3)$ accuracy. The extension can be easily done either continuously or discretely.  For a Poisson equation, we just need to extend one $f$ value
for both 2D and 3D problems.  If one prefers not to  extend the source term, then there is not enough degree of freedom for fourth order compact schemes. Our HOC methods can ensure a super-third order  convergence   that are accurate to $O(h^4)$ for all fourth order polynomials except for the $x^4$ and  $y^4$ terms, which is referred  as  the super-third order  methods in this paper.

We think the proposed new strategy is quite versatile. Just a few lines of modifications to the computer codes are needed for the different  problems mentioned above. We also think that  most of HOC schemes in the literature can be included in the framework of our method  from the definition of HOC schemes even though the derivations may be  different. The designed HOC schemes are recommended for Poisson/Helmholtz equations, diffusion and advection equations, anisotropic  PDEs with constant coefficients since the finite difference coefficients just need to be computed once or twice.  The proposed methods can be applied to PDEs with variable coefficients, but the computational costs may be too high to be practical. For variable coefficients problems, a Richardson extrapolation approach \cite{TONG2021107413} may be the most economical way to get a fourth order method, although it would not be compact anymore.

 The developed new methods can be and have been applied to  Helmholtz type equations  ($K>0$). When $K\le 0$, then the proposed methods always work with strict error bounds; \ if $K> 0$, then the proposed methods  always work if $h$ is smaller enough. 
The computed solutions have the designed order of convergence but the error constant is  promotional to $K$. 
If  $h$ is not smaller enough and $K$ is large, the developed methods work fine as long as the coefficient matrix $A_h$ is not singular; that is, $K$ is not an eigenvalue of $A_h$, for which the probability is one. We do not have a uniform error bound in this case though and the error depends on $\|(A_h)^{-1}\|$.  

The rest of the paper is organized as follows. In the next section, we explain the algorithm for Poisson/Helmholtz equations with a Neumann/Robin boundary condition along part of a boundary.  In Section~\ref{sec:super-3rd}, we explain the super-third algorithm and present numerical test results and comparisons.  In Section~\ref{sec:diff-advec}, we discuss the fourth order compact scheme for   diffusion and advection equations  with constant coefficients.  In Section~\ref{sec:anis}, we discuss the HOC schemes for anisotropic diffusion and advection equations for Dirichlet (fourth-order) and flux (super-third) boundary conditions.
In  Section~\ref{sec:3D}, we discuss the fourth order compact scheme for flux
boundary conditions in three dimensions followed with a numerical example.
We conclude in the last section.

\section{Constructing 4th-order compact  schemes for Poisson/Helmholtz equations with a flux BC using f-extension}

We first construct a fourth order compact scheme for a  Helmholtz (including Poisson when $K=0$)  equation on a rectangular domain in 2D. Without loss  of generality and for the convenience of presentation, we assume that the domain $\cal{R}$ is a  square  $[x_l, \,x_r]\times [y_l, \,y_r]$.
 We use a uniform  mesh
\eqm
  x_i =  x_l + ih, \quad i=0,1,\cdots, N;  \qquad y_j=  y_l+jh, \quad j=0,1,\cdots, N.
\enm
\begin{figure}[phbt]
\centerline{  \includegraphics[width=0.65\textwidth]{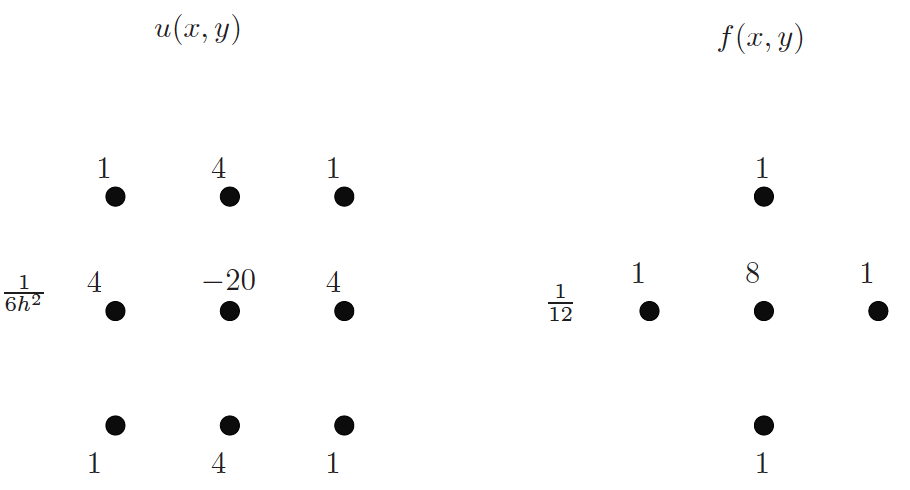} }
\caption{The coefficients of the finite difference scheme $U_{ij}$ using the compact
nine-point stencil (left diagram) and the linear combination of $f(x_i,y_j)$ \ (right diagram).}
\label{9lap}
\end{figure}
\noindent
In an interior grid point,  the classical compact nine-point finite difference stencil, the finite difference coefficients and the right hand sides are illustrated  in Figure~\ref{9lap}, see for example, \cite{strikwerda,li:book-FDFEM}.
The fourth order compact scheme can be written as
\eqml{compact2d}
\lp  L_h + K M_h\rp U_{ij}   = M_h f_{ij},
\enml
at {\em an interior grid point}  $(x_i, y_j)$, where $L_h$ is the discrete nine-point  Laplacian whose coefficients at four corners are $\frac{1}{6 h^2}$, at
east-north-south-west grid points are $\frac{4}{6 h^2}$, and at the center is $-\frac{20}{6 h^2}$; $M_h$ is an averaging  operator whose coefficients at east-north-south-west grid points are $\frac{1}{12}$, and at the center is $\frac{8}{12}$, that is,
\eqm
  M_h f_{ij} = \frac{1}{12} \left ( \frac{\null}{\null} f_{i-1,j} + f_{i+1,j}+f_{i,j-1} + f_{i,j+1} + 8 f_{i,j}\right),
\enm
where $f_{ij} = f(x_i,y_j)$ and so on, see Figure~\ref{9lap} for an illustration.
Assume that the coefficient matrix of the finite difference equation is $A_h$.  Then, $-A_h$ is an M-matrix if $K\le 0$ and  $4/ (6h^2) + K/12\ge 0$, {\em i.e.} $1/h^2 + K/8\ge 0$.
%
There are several definitions of an M-matrix in the literature. 
 We quote the definition from \cite{hackbusch1994iterative} below.
\begin{definition}
  Let $A$ be  $n\times n$ real matrix. Then matrix $A$ is called an M-matrixx if
  \begin{align}
   & a_{ii} > 0  \textrm{ for all } 1 \leq i\leq  n, \label{signa} \\
   & a_{ij} \leq 0 \textrm{ for all } i\neq j, 1 \leq i,j \leq  n, \label{signb} \\
   & A \; \textrm{is regular and } A^{-1}\geq 0. \label{signc}
  \end{align}
\end{definition}

The sign conditions (\ref{signa}) and (\ref{signb}) are easy to check, but not so for the condition (\ref{signc}).   
We  combine Theorem~6.4.10 and Lemma 6.4.12 from \cite{hackbusch1994iterative}  to have the following lemma to state an equivalent definition of an M-matrix.

\begin{lemma} \label{lem1}
  Let the matrix $A\in R^{n\times n}$ be irreducibly weakly diagonally dominant with at least one row being strictly. If the sign conditions (\ref{signa}) and (\ref{signb})
  are satisfied, then $A$ is an M-matrix.  Furthermore, if $A=A^T$, then $A$ is also symmetric positive definite.
\end{lemma}

From the definitions of $L_h$, $M_h$ in (\ref{compact2d}) and Lemma \ref{lem1}, we know that the coefficient matrix ($-A_h$) of the finite difference equations is an M-matrix  
if $K\le 0$  and $1/ h^2 + K/8\ge 0$\footnote{This is because the $Ku$ term is treated as a source. Note that if the condition is violated, that means  $|K|> 8/h^2,\, K<0$, then the coefficient matrix is actually more diagonally dominant since the contribution to the diagonal is $8K/12$ while to each of the four off-diagonals is $K/12$.   The coefficient matrix is invertible and the norm $\|A_h^{-1}\|_2$ will be smaller than that corresponding to $K=0$ although the sign property will be violated. In our new methods, we do not treat $Ku$ as a source term and just need the usual $K\le 0$ condition for an M-matrix condition.}.
%
%
The local truncation error at {\em an interior grid point}  $(x_i, y_j)$  is defined as
\eqm
  T_{ij}^h = \lp L_h + K M_h \rp u(x_i,y_j)    -M_h f(x_i,y_j),
\enm
which is of $O(h^4)$ if $u\in C^6(\cal{R})$. Thus,  the global error is bounded by  $\|E_h\|_{\infty} =\|u(x_i,y_j)-U_{ij}\|_{\infty} \le C h^4$   if  a Dirichlet boundary condition is specified. 

Now we  consider a Robin boundary condition $\left. ({\frac{\partial u}{\partial n}}  +  \sigma u) \right |_{x=x_0} = g(y)$  to explain our idea and the new algorithm. Note that ${\bf n} = (-1,0)$ and $\frac{\partial u}{\partial n}  = -\frac{\partial u}{\partial x}$ at the boundary $x=x_0$. We assume that we know the source term $f(x,y)$ at $x=x_0-h$, or we can use a quadratic extension that is third order accurate . In the next section, we will explain a super-third order compact algorithm that does not need to use $f(x,y)$ at $x=x_0-h$.

A compact FD scheme at a boundary grid point $(x_0,y_j)$, $j=1,2,\cdots,N-1$ can be written as
\eqml{4th-hoc}
&  \dsp  \sum_{i_k=0}^{1} \sum_{j_k=-1}^{1} \!\!  \alpha_{i_k,j_k} U_{i+i_k,j+j_k}   =  \sum_{i_k=-1}^{1} \sum_{j_k=-1}^{1}  \!\!\beta_{i_k,j_k} \, f(x_{i+i_k},y_{j+j_k})  + \sum_{j_k=-1}^{1} \!\! {\gamma_{j_k}} \, g(y_{j+j_k})  , \\
  & \dsp  \sum_{i_k,j_k=-1}^1 \!\!\! \beta_{i_k,j_k} = 1,
\enml
where $i=0$ throughout this section, $\alpha_{i_k,j_k}$ ($i_k=0,1$, $j_k=-1,0,1$), $\beta_{i_k,j_k} $ ($i_k=-1,0,1$, $j_k=-1,0,1)$, and ${\gamma_{j_k}}$($j_k=-1,0,1$) are undetermined coefficients. We leave the index $i$ in the formulas  so that the derivation can be applied to interior grid points as well for other problems discussed in the later sections. Thus, we have $6$ coefficients for $U_{ij}$,  $9$ coefficients for $f_{ij}$, {and 3 coefficients for $g(y_j)$.}  The reason to have the constraint $\dsp \sum_{i_k,j_k=-1}^1 \beta_{i_k,j_k} = 1$ will be seen soon, which is analogues to the fourth order scheme at interior grid points. Apparently, almost all  existing high order compact schemes have the form above.

\ignore{
\begin{remark}
 We think all the existing high order compact schemes should have the form above. The derivation of our method can lead to other high order compact schemes except the constraint  in the second identity may be different.
\end{remark} }

Denote the `local truncation error' for the scheme at  $(x_i,y_j)$ as
\eqmno
  T_{ij}^h = \sum_{i_k=0}^{1} \sum_{j_k=-1}^{1}  \!\!  \alpha_{i_k,j_k} u(x_{i+i_k},y_{j+j_k} ) -  \sum_{i_k=-1}^{1}  \sum_{j_k=-1}^{1} \!\!  \beta_{i_k,j_k} f(x_{i+i_k},y_{j+j_k})  - \sum_{j_k=-1}^1 \! \!\!\gamma_{j_k}  {g(y_{j+j_k})} .
\enmno
We want to determine the coefficients such that the local truncation errors are zero or $O(h^4)$ for any fourth order polynomials. To derive the system of equations for the coefficients, we expand $u(x_{i+i_k},y_{j+j_k} )$, $f(x_{i+i_k},y_{j+j_k})$ at $(x_i,y_j)$, and ${g(y_{j+j_k})}$ at $y_j$. For the $u(x_{i+i_k},y_{j+j_k} )$ terms, we apply the Taylor expansion at $(x_i,y_j)$ up to all  fourth order partial derivatives,
\eqmno
  u(x_{i+i_k},y_{j+j_k} ) =  \sum_{ 0 \leq k_1+k_2 \leq 4} \!\!  \left.  \frac{1}{k_1!\, k_2! }\, \frac{\partial^{k_1+k_2}  u}{\partial x^{k_1} \partial y^{k_2} } \right |_{(x_i,y_j)}
  \!\! \!\! h_{i_k}^{k_1} h_{j_k}^{k_2} + O( h^5),
\enmno
where $h_{i_k} = i_k h$ and $h_{j_k}= j_k h$.

 For the $f(x_{i+i_k},y_{j+j_k} )$ terms, we apply the Taylor expansion at $(x_i,y_j)$ up to  all   second order partial derivatives,
\eqm
  f(x_{i+i_k},y_{j+j_k} ) = f_0 + f_{x,0} \,h_{i_k} + f_{y,0} \,h_{j_k} + f_{xx,0} \, \frac{h_{i_k}^2}{2}  +  f_{xy,0} \, h_{i_k} h_{j_k} + f_{yy,0}\,  \frac{h_{j_k}^2}{2}  + O(h^3),
  \enm
where $f_0=f(x_i,y_j)$,  $f_{x,0}= \frac{\partial f}{\partial x}(x_i,y_j)$ and so on.  By differentiating the Helmholtz equation, we get the following high order PDE relations
\eqml{fxxetc}
 & \dsp f_x = u_{xxx} + u_{yyx} + K u_x, \qquad f_y = u_{xxy} + u_{yyy} + K u_y, \\ \eqsp
  & \dsp f_{xx} = u_{xxxx} + u_{yyxx} + K u_{xx} , \qquad  f_{xy} = u_{xxxy} + u_{yyyx} + K u_{xy}  , \\ \eqsp
  &  f_{yy} = u_{xxyy} + u_{yyyy} + K  u_{yy} .
\enml

For the $g(y_{j+j_k})$ terms, we apply the Taylor expansion at $y_j$ with respect to $y$ up to third  order  derivatives
\eqm
g(y_{j+j_k}) = g(y_{j}) + g'(y_{j})h_{j_k} + g''(y_{j})\frac{h_{j_k}^2}{2} + g'''(y_{j})\frac{h_{j_k}^3}{6} + O(h^4).
\enm
Note that,
\eqm
  \left .  \lp \dsp \frac{\partial u}{\partial n}  +  \sigma  u (\mathbf{x})  \rp \right |_{x=0} =  \left .  \lp \dsp  - \frac{\partial u}{\partial x}   +  \sigma  u (\mathbf{x})  \rp \right |_{x=0}=  g(y).
\enm
Thus, we have more relations from the boundary condition needed for the HOC scheme,
\eqml{mix_u}
 && \dsp -u_x(0,y) + \sigma  u (0,y) = g(y), \qquad  \qquad  -u_{xy} (0,y) + \sigma  u_y(0,y) = g'(y),  \\ \eqsp
 && \dsp   - u_{xyy} (0,y) + \sigma  u_{yy} (0,y) = g''(y),
 \qquad  - u_{xyyy} (0,y) + \sigma  u_{yyy} (0,y) = g'''(y).
\enml
These relations are utilized in deriving the fourth order compact scheme.

After the expansions of all involved terms, we can write  the local truncation error as
\begin{equation}\label{local-error}
  T_{ij}^h=\sum_{ 0 \leq k_1+k_2 \leq 4}  \!\!\!\!\! L_{k_1,k_2} \, \left.   \frac{\partial^{k_1+k_2}  u}{\partial x^{k_1} \partial y^{k_2} } \right |_{(x_i,y_j)}
  + O\left(\| {\bfalf} \|_{\infty} h^5 + \| {\bfbeta} \|_{\infty} h^3 +\| {\bfgamma} \|_{\infty} {h^4}\right ),
\end{equation}
where $L_{k_1,k_2}$ are the results after we collect terms and will be seen in the linear system of equations for the coefficients soon, $\| {\bfalf} \|_{\infty} = \max\limits_{0\leq i_k \leq 1,-1\leq j_k \leq 1}\left\{|\alpha_{i_k,j_k}|\right\}$ and so on.
We want the local truncation error to be zeros or $O(h^4)$ for all fourth order polynomials,  or  $x^{k_1} y^{k_2}$, $0 \leq k_1+k_2 \leq 4$. The finite difference equation should approximate the Poisson equation well for which $f(x,y)$ is an $O(1)$ quantity in general. Thus, as in the standard 9-points fourth order compact scheme, we impose the constraint of $\dsp\sum_{i_k,j_k=-1}^1 \!\! \beta_{i_k,j_k} = 1$ otherwise the minimum of  $|T_{ij}^h|$  are zero with all zero coefficients. 
In this way, by matching the terms of the coefficients of $u, u_x, u_y, \cdots, u_{xxxx}, u_{xxxy}, \cdots, u_{yyyy}$, we have $15$ linear equations with one constraint. The first six equations (required for all quadratic polynomials) are,
\eqml{poisson-coeA}
&& \dsp \sum_{i_k=0}^{1} \sum_{j_k=-1}^{1} \alpha_{i_k,j_k}  - K -   \sigma  \sum_{j_k=-1}^1 \gamma_{j_k}   =0 \\ 
&&\dsp  \sum_{i_k=0}^{1} \sum_{j_k=-1}^{1} \alpha_{i_k,j_k}  h_{i_k}  -  K \!\!  \sum_{i_k=-1}^{1} \sum_{j_k=-1}^{1}  \beta_{i_k,j_k} h_{i_k} {+} \sum_{j_k=-1}^1 \!\! \gamma_{j_k}   =0 \\ 
&&\dsp  \sum_{i_k=0}^{1} \sum_{j_k=-1}^{1} \alpha_{i_k,j_k}  h_{j_k}  -  K \!\!  \sum_{i_k=-1}^{1} \sum_{j_k=-1}^{1}  \beta_{i_k,j_k} h_{j_k}  -  \sigma  \!\! \sum_{j_k=-1}^1 \!\! \gamma_{j_k}  h_{j_k}  =0  \\ 
&& \dsp \sum_{i_k=0}^{1} \sum_{j_k=-1}^{1} \alpha_{i_k,j_k}  \frac{h_{i_k}^2 } {2} -  \sum_{i_k=-1}^{1} \sum_{j_k=-1}^{1}  \beta_{i_k,j_k}
-  K \!\!  \sum_{i_k=-1}^{1} \sum_{j_k=-1}^{1}  \beta_{i_k,j_k}  \frac{h_{i_k}^2 } {2} = 0 \\ 
&& \dsp \sum_{i_k=0}^{1} \sum_{j_k=-1}^{1} \alpha_{i_k,j_k}  \frac{h_{j_k}^2 } {2} -  \sum_{i_k=-1}^{1} \sum_{j_k=-1}^{1}  \beta_{i_k,j_k}
-  K\!\!  \sum_{i_k=-1}^{1} \sum_{j_k=-1}^{1}  \beta_{i_k,j_k}  \frac{h_{j_k}^2 } {2}  -  \sigma  \!\! \sum_{j_k=-1}^1 \!\! \gamma_{j_k}  \frac{h_{j_k}^2 } {2}   = 0 \\  
 && \dsp \sum_{i_k=0}^{1} \sum_{j_k=-1}^{1} \alpha_{i_k,j_k} h_{i_k} h_{j_k}
   -  K \!\!  \sum_{i_k=-1}^{1} \sum_{j_k=-1}^{1}  \beta_{i_k,j_k} h_{i_k}  h_{j_k}
 {+} \sum_{j_k=-1}^1 \!\! \gamma_{j_k}  h_{j_k} = 0.
 \enml
 The next four equations (required for cubic polynomials) are
\eqml{poisson-coeB}
&& \dsp \sum_{i_k=0}^{1} \sum_{j_k=-1}^{1} \alpha_{i_k,j_k}  \frac{h_{i_k}^3 } {3!} -  \sum_{i_k=-1}^{1} \sum_{j_k=-1}^{1}  \beta_{i_k,j_k} h_{i_k} = 0 \\ \eqsp
&& \dsp \sum_{i_k=0}^{1} \sum_{j_k=-1}^{1} \alpha_{i_k,j_k}  \frac{h_{i_k}^2  h_{j_k}} {2} -  \sum_{i_k=-1}^{1} \sum_{j_k=-1}^{1}  \beta_{i_k,j_k} h_{j_k} = 0 \\ \eqsp
  && \dsp \sum_{i_k=0}^{1} \sum_{j_k=-1}^{1} \alpha_{i_k,j_k}  \frac{h_{i_k}  h_{j_k}^2} {2} -  \sum_{i_k=-1}^{1} \sum_{j_k=-1}^{1}  \beta_{i_k,j_k} h_{i_k}  {+} \sum_{j_k=-1}^1 \gamma_{j_k}    \frac{h_{j_k}^2 } {2} = 0  \\ \eqsp
&&\dsp  \sum_{i_k=0}^{1} \sum_{j_k=-1}^{1} \alpha_{i_k,j_k}  \frac{  h_{j_k}^3} {3!} -  \sum_{i_k=-1}^{1} \sum_{j_k=-1}^{1}   \!\!  \beta_{i_k,j_k} h_{j_k}  -  \sigma \sum_{j_k=-1}^1  \!\!  \gamma_{j_k}   \frac{h_{j_k}^3 } {3!}  = 0.
\enml

The next  five equations (required for quartic polynomials) are
\eqmno 
&& \dsp \sum_{i_k=0}^{1} \sum_{j_k=-1}^{1} \alpha_{i_k,j_k}  \frac{h_{i_k}^4 } {4!} -  \sum_{i_k=-1}^{1} \sum_{j_k=-1}^{1}  \beta_{i_k,j_k}  \frac{h_{i_k}^2  } {2}  = 0 \\ 
&& \dsp \sum_{i_k=0}^{1} \sum_{j_k=-1}^{1} \alpha_{i_k,j_k}  \frac{h_{i_k}^3  h_{j_k}} {3!} -  \sum_{i_k=-1}^{1} \sum_{j_k=-1}^{1}  \beta_{i_k,j_k} h_{i_k} h_{j_k} = 0 
 \enmno 
\eqml{poisson-coeC}
   && \dsp \sum_{i_k=0}^{1} \sum_{j_k=-1}^{1} \alpha_{i_k,j_k}  \frac{6 h_{i_k}^2 h_{j_k}^2} {4!} -  \sum_{i_k=-1}^{1} \sum_{j_k=-1}^{1}  \beta_{i_k,j_k}  \left (  \frac{h_{i_k}^2 } {2} +  \frac{h_{j_k}^2 } {2} \right ) = 0 \\ 
&& \dsp \sum_{i_k=0}^{1} \sum_{j_k=-1}^{1} \alpha_{i_k,j_k}  \frac{h_{i_k}  h_{j_k}^3} {3!} -  \sum_{i_k=-1}^{1} \sum_{j_k=-1}^{1}  \beta_{i_k,j_k} h_{i_k} h_{j_k}  {+} \sum_{j_k=-1}^1 \gamma_{j_k}   \frac{h_{j_k}^3 } {6} = 0 \\ 
&& \dsp \sum_{i_k=0}^{1} \sum_{j_k=-1}^{1} \alpha_{i_k,j_k}  \frac{  h_{j_k}^4} {4!} -  \sum_{i_k=-1}^{1} \sum_{j_k=-1}^{1}  \beta_{i_k,j_k}  \frac{h_{j_k}^2 } {2} = 0  .
 \enml

\begin{figure}
\centering
\!\!\! \subfigure[$U_{ij}$]{\includegraphics[width=0.22\textwidth]{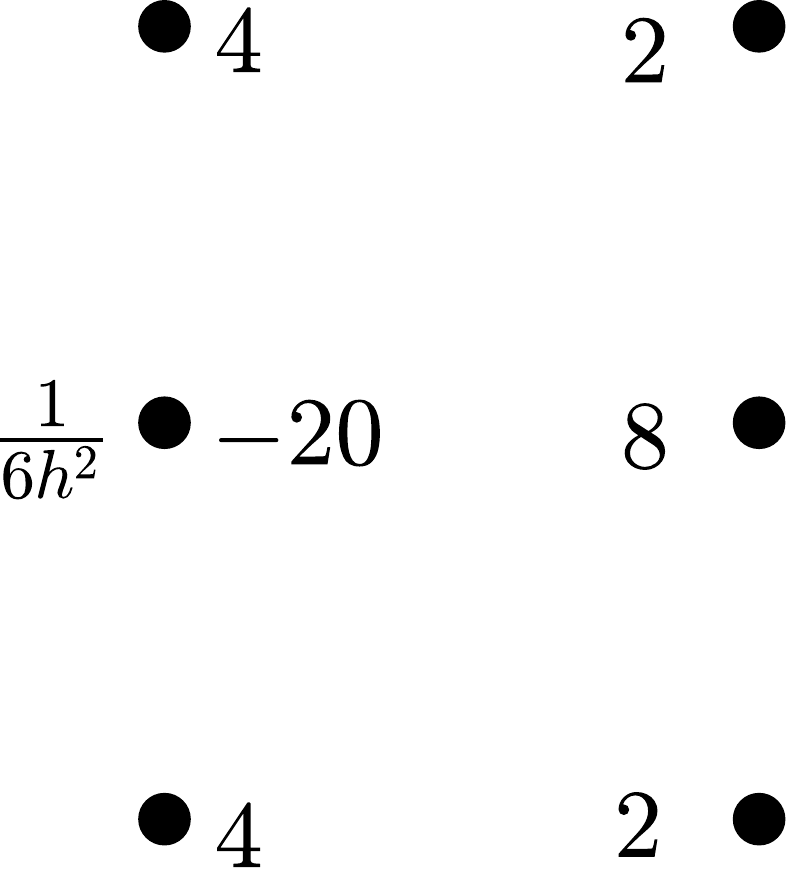}}
$\null \; \qquad\; \null \quad$$\null \quad$
\subfigure[$f _{ij}$]{\includegraphics[width=0.26\textwidth]{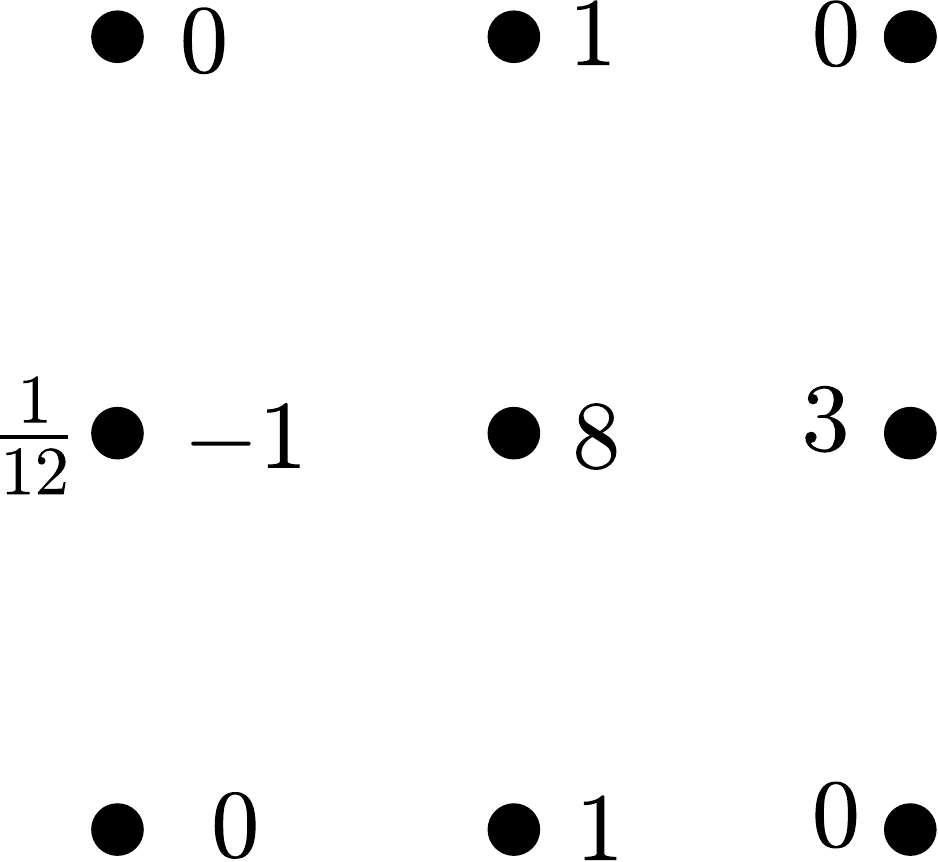}} $\null \; \qquad\; \null \quad$
$\null \quad$
\subfigure[${-g_j}$]{\includegraphics[width=0.044\textwidth]{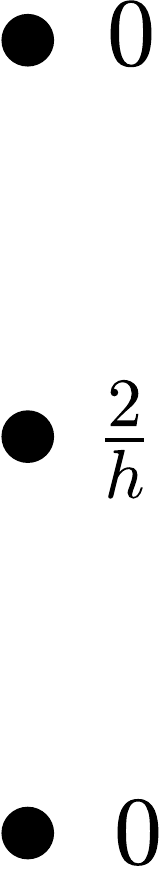}}
\caption{A set of coefficients of the fourth order compact  scheme for  Poisson equations with a Neumann BC given at the left boundary $x=x_l$. The left diagram is for that of $U_{ij}$; the middle is for $f_{ij}$;  and the right is for $g_{j}$.  For a Robin BC, we need to add $-\dsp {2 \sigma}/{h}$ to the coefficient $\alpha_{0 j}$ for $U_{0j}$, which is on a diagonal of the coefficient matrix of the FD equation.}
\label{hoc_diag_Poisson}
\end{figure}

There are $16$ equations and $18$ unknowns, which is an under-determined system of equations. In general there are infinity number of solutions.
We have found out {{\em a set of coefficients analytically}} \ for Poisson equations with a Neumann BC, 
given at the left boundary $x=x_l$, see Figure~\ref{hoc_diag_Poisson},
which lists a set of  coefficients $\alpha_{i_k,j_k}$ corresponding to the finite difference coefficients $U_{ij}$; $\beta_{i_k,j_k}$ corresponding to the combination coefficients of $f$; and $\gamma_{j_k}$ corresponding to the combination coefficients of $g$.

For a Poisson equation with a Robin BC, we need to add $-\dsp {2 \sigma}/{h}$ to the coefficient $\alpha_{0,j}$ for $U_{0j}$, which is on a diagonal of the coefficient matrix of the FD equations. 
Note also that,  similar to the ghost point method,  those coefficients at $(x_0+h,y_j)$ are doubled  compared with those at an interior grid point, see the corresponding coefficients in Fig.~\ref{9lap} and Fig~\ref{hoc_diag_Poisson}, which  can be  regarded as the reflection, shifting the coefficient at $(x_0-h,y_j)$ to $(x_0+h,y_j)$ with an appropriate adjustment for $f(x_0-h,y_j)$ and $f(x_0+h,y_j)$ compared with that of an interior grid point.
It is worth mentioning that, for a Poisson equation and other PDEs with constant coefficients with a flux BC,
the derivation of the scheme is independent of the index $j$. Thus,  we can simply use the grid points $(0, \pm h)$,  $(h, \pm h)$, $(h,0)$,  $(0,0)$, the six particular grid points to derive the coefficients.
Note also that there are other analytic sets of solutions as well.


For Helmholtz equations with a Robin boundary condition, we have also found a set of analytic coefficients  below obtained from the  Maple symbolic package:
\eqm \label{Coeff_H_A}
& U_{i,j}:  \;  \dsp  \frac{1}{\lambda h^2}\left[
               \begin{array}{cc}
                  8 & 4 \\
                  4\sigma K h^3 + \lambda K h^2 - 24 \sigma h - 40 \null \quad & 16 \\
                  8 & 4 \\
               \end{array}
             \right], \\ \eqsp
& f_{i,j}: \; \dsp \frac{1}{\lambda}\left[
               \begin{array}{ccc}
                 0 & 1 & 0 \\
                 -1 & 8-Kh^2 & 3 \\
                 0 & 1 & 0 \\
               \end{array}
             \right],   \qquad
g_{j}: \;   \dsp \frac{1}{\lambda h}\left[
               \begin{array}{c}
                  0  \\
                  4Kh^2-24  \\
                  0  \\
               \end{array}
             \right],  \label{Coeff_H_B}
\enm
with $\lambda = 12-Kh^2$, corresponding to the finite difference coefficients for $U_{ij}$, the combination coefficients for $f_{ij}$, and  $g_{j}$ of the right hand side of the flux boundary condition. Again we see the perfect symmetry in the coefficients. When $K=0$, the finite difference scheme is the same as that for Poisson equations; and when $K\le 0$, the coefficient matrix of the finite difference equations is an M-matrix and the convergence is guaranteed.
While there are infinite number of solution sets to the linear system of equations for the coefficients, we show that the system of equations for the coefficients are solvable.


Next, we discuss the convergence of the fourth order compact scheme for Poisson or generalized Helmholtz equations ($K< 0$) as summarized in the following theorem.
The sign property, $\alpha_{00} <0$ and otherwise $\alpha_{i,j} \ge 0$ is the key for the convergence proof.
In general, there is no guarantee for the convergence of Helmholtz equations ($K>0$) since there is no control of the coefficient matrix of the FD equations. For instance, when $K$ is in a neighborhood of an eigenvalue of the boundary value problem, then the  coefficient matrix  is close to be singular.

\begin{theorem} \label{thm-A}
Let $U_{ij}$ be the finite difference solution obtained from the fourth order scheme with the set of the coefficients
given by  (\ref{Coeff_H_A})-\eqn{Coeff_H_B}, which include the special situation listed in Figure~\ref{hoc_diag_Poisson}.
 Assume that the solution to the Poisson/Helmholtz
equation with non-empty ${\cal{R}}_1$(Dirichlet) and  ${\cal{R}}_2$(Robin) in the boundary condition~\eqn{two-BC} is $u(x,y)$. Then, the algorithm is {\em exact}  if the solution is  a fourth  order polynomial if $K=0$ and $\sigma=0$.  For a general solution
  $u(x,y)\in C^6(\cal{R})$,   we have the following error estimate assuming that $K\le 0$ and $\sigma \ge 0$,
\eqm
  \left \| u(x_i,y_j)- U_{ij} \right \|_{\infty}  \le C \, max\{ \|D^6 u\|_{\infty, {\cal{R}}},  \, \|D^5 u\|_{\infty, \partial {\cal{R}}}\} h^4,
\enm
\end{theorem}
 where $D^6 u$ means all possible sixth order partial derivatives of $u(x,y)$ and so on. 

\vthin

{\bf Proof:} For the fixed set of the coefficients, we know that the scheme is exact for all fourth order polynomials at all (interior and boundary) grid points  from the design of the algorithm if $K=0$ and $\sigma=0$.

For a general solutions  $u(x,y)\in C^6(\cal{R})$, we know that the local truncation errors are of $O( \|D^6 u\|_{\infty, {\cal{R}}} \,h^4)$ at interior grid points and are of $O(\|D^5 u\|_{\infty, \partial {\cal{R}}} \,h^3)$ at  boundary grid points $(x_0, y_j)$ from  (\ref{local-error}) and the analytic expressions of the coefficients, which implies that $ {\bfalf}_{\infty}= O(\frac{1}{h^2})$, $ {\bfbeta}_{\infty}= O(1)$ and $ {\bfgamma}_{\infty}= O(\frac{1}{h})$.

If $K\le 0$ and $\sigma \ge 0$, from (\ref{Coeff_H_A}) and Lemma \ref{lem1}  we know that the coefficient matrix ($-A_h$) of the finite difference equations is an M-matrix. Thus, from  Theorem~6.1 and Theorem~6.2 of Morton~\&~Mayer's book \cite{morton-mayers}, which is also valid for part of Neumann boundary condition as long as at least  one point has a Dirichlet boundary condition prescribed as stated in the book,  we conclude that $\left \| u(x_i,y_j)- U_{ij} \right \|_{\infty}  \le C h^4$. Note that, for part of Neumann or Robin boundary conditions, the related boundary grid points are regarded as an `interior points' in the proof of the maximum principle, which implies a slightly larger error constant. 
\hfill $\square$

The  fourth order compact scheme is not exact for Helmholtz equations because  we only  utilized up to all second order
partial derivatives for the $Ku$ term.  Note also that for a pure Neumann boundary condition along the entire boundary, we need to assume that the solution exists. Then, we can specify the solution at one point at the boundary to make the solution unique and  apply the theorem. For other well-posed situations such as $K< 0$, or part of Robin BC's without a Dirichlet BC on the boundary, we  believe that the method is still fourth order accurate, but the  proof is an open challenge.

\subsection{Numerical experiments of the fourth order compact scheme for Poisson and Helmholtz equations}

We have tested the proposed fourth order compact scheme for polynomials $P_k(x,y)$, $k\le 4$, the computed solutions are accurate to  $\epsilon \,\mbox{cond}(A_h)$ when $K=\sigma=0$, where $\epsilon\sim 10^{-16}$  is the machine precision and  $A_h$ is the coefficient matrix of the   finite difference equations.
Next we  test  two constructed examples with genuine non-linear solutions, in which  one  is a relatively smooth, and the other  can be   oscillatory.

\begin{example} \label{ex1a} An example with a smooth solution.
\eqml{ex1eq}
 &&  \dsp u(x,y) = e^{-x} \sin (\pi y), \quad  \quad  (x,y) \in (0,\;1)^2, \\ \eqsp
  && \dsp f(x,y) = e^{-x}  \sin(\pi y) (1-\pi^2 ),  \quad  (x,y) \in (0,\;1)^2, \\ \eqsp
 && \dsp  \left. \lp \frac{\null}{\null}  \frac{\partial u}{\partial n} +  \sigma u \rp  \right |_{x=0} = (1+\sigma )  \sin ( \pi y) , \qquad y\in  (0,\;1).
 \enml
\end{example}
In this example, the Robin BC is a non-zero function of $y$. The solution and the source term are relatively smooth.


\begin{example} \label{ex2}
An example with an oscillatory solution.
\eqml{ex1}
 &&  \dsp u(x,y) = \sin(k_1 x) \cos(k_2 y), \quad \quad  (x,y) \in (0,\;1)^2, \\ \eqsp
 &&  \dsp f(x,y) =  - \left (k_1^2 + k_2^2 \right )   \sin(k_1 x) \cos(k_2 y),  \quad  (x,y) \in (0,\;1)^2, \\ \eqsp
 && \dsp  \left. \lp \frac{\null}{\null}  \frac{\partial u}{\partial n} + \sigma u \rp  \right |_{x=0} =   -k_1  \cos(k_2 y)
 +  \sigma \sin(k_1 x) \cos(k_2 y), \qquad y\in  (0,\;1)  .
 \enml
\end{example}
In this example, we can choose $k_1$ and $k_2$ to make the solution more oscillatory. For a typical test, we choose $k_1=5$ and $k_2=50$, so $f(x,y)\sim 2525$. This is a relatively tough problem to compute for large $k_1$ or $k_2$ since $\frac{\partial^5 u}{\partial y^5} \sim  3.125\times 10^{8}$. The mesh needs to be fine enough to resolve the solution.

We assume Dirichlet boundary conditions on other parts of the boundary  from  the exact solution.
In  Table~\ref{tab:4thA}, we show some  experimental results. The top table lists results for Example~\ref{ex1a}.  The second-third columns are the results for the Poisson equation with a Neumann BC, while the fourth-fifth columns list the results for  the Helmholtz  equation with $K=2000$ and a Robin boundary condition with $\sigma=-20<0$ to test our method for an extreme case.
In the table, $N$ is the number of grid lines in one coordinate direction  so $h=1/N$; and the order is the computed convergence order using two consecutive errors,
\eqm
  order =\frac{ \log ( \|E_N\|_{\infty} /\log \|E_{2N} \|_{\infty} )}{\log 2} .
\enm
We see clearly fourth order convergence. From $N=256$ to $N=512$, we observe better than expected convergence order for which we think it just a coincidence. The results  using $N=510$ or $N=513$ are in  line of a fourth order method. Also note that  around $N=512$, the mesh $h$  is close to the best possible  before the round-off errors  become dominant to  ruin the convergence  if $h$ decreases further.

\begin{table}[htbp]
\caption{Grid refinement analysis  of the fourth order compact scheme.   The top table lists the results for Example~\ref{ex1a}  with a Neumann BC and $K=0$ in the column 2-3, and a Robin BC ($\sigma=-20$) in the columns 4-5. The bottom table lists the results for Example~\ref{ex2}  for  the
Helmholtz equation with $K=2000$, $k_1=5$, $k_2=50$, and a Neumann BC in the columns 2-3, and a Robin BC ($\sigma=-20$) in the columns 4-5. In both cases, we see clearly fourth order convergence.} \label{tab:4thA}

\vthin

\begin{center}
\begin{tabular}{|c| c c | c  c|}
  \hline
       $N$   & $\|E\|_{\infty} \,   Neumann  $  & $order $    & $\|E\|_{\infty} \; \;  {\color{blue} Robin} $  & $order $  \\       \hline	
  $16$  &  2.2943\,e-05  &  $ $ &  1.6933e-05            & \\
  $32$  &  1.4127\,e-06  & $4.0215 $ & 1.0589e-06  &  3.9992 \\
  $64$  &  8.7602\,e-08  & $ 4.0114 $ & 6.6187e-08  &  3.9999 \\
 $128$  & 5.4524\,e-09   & $4.0060$ & 4.1360e-09 &  4.0002\\
 $256$  & 3.3800\,e-10  & $4.0118$  & 2.6803e-10  &  3.9478 \\
 $512$  &  2.1125\,e-11  & $ 4.7221$ & 1.3337e-11  &  4.3289 \\
      \hline
\end{tabular}

\vthin

\begin{tabular}{|c| c c | c c|}
  \hline
       $N$   & $\|E\|_{\infty} \;  Neumann $  & $order$   & $\|E\|_{\infty} \; \; {\color{blue} Robin} $  & $order $   \\       \hline	
  $16$  &  1.0920\,e-00  & &  2.8336  & \\
  $32$  &  4.7754\,e-02  & $  4.5152$ & 1.5326e-01  & 4.2086 \\
  $64$  &  2.7619\,e-03  & $ 4.1119 $  & 7.4075e-03  &  4.3709 \\
 $128$  & 1.6941\,e-04   & $4.0271$ & 4.2780e-04  &  4.1140\\
 $256$  & 1.0539\,e-05  & $ 4.0067$  & 2.6284e-05  & 4.0247 \\
 $512$  &  6.5794\,e-07  & $  4.0016$ & 1.6326e-06  &  4.0089 \\
      \hline
\end{tabular}

\end{center}
\end{table}

In the bottom table of Table~\ref{tab:4thA}, we show the grid refinement analysis  for the Helmholtz equation  using  Example~\ref{ex2} with $K=2000$, $k_1=5$, $k_2=50$.  The second-third columns are the results for the  Neumann BC, while the fourth-fifth columns list the results for  a Robin boundary condition with $\sigma=-20$.
We see errors are larger compared with that for Example~\ref{ex1a}  due to the nature of the oscillatory solutions and large amplitudes of high order partial derivatives. The mesh  needs to  be fine enough to resolve the solution. Still, we see clearly fourth  order convergence. Also, the coefficient $\sigma$ has little effect on the convergence unless it is very large.

\begin{figure}[phbt]
\begin{minipage}[t]{3.0in} (a)

\includegraphics[width=1.0\textwidth]{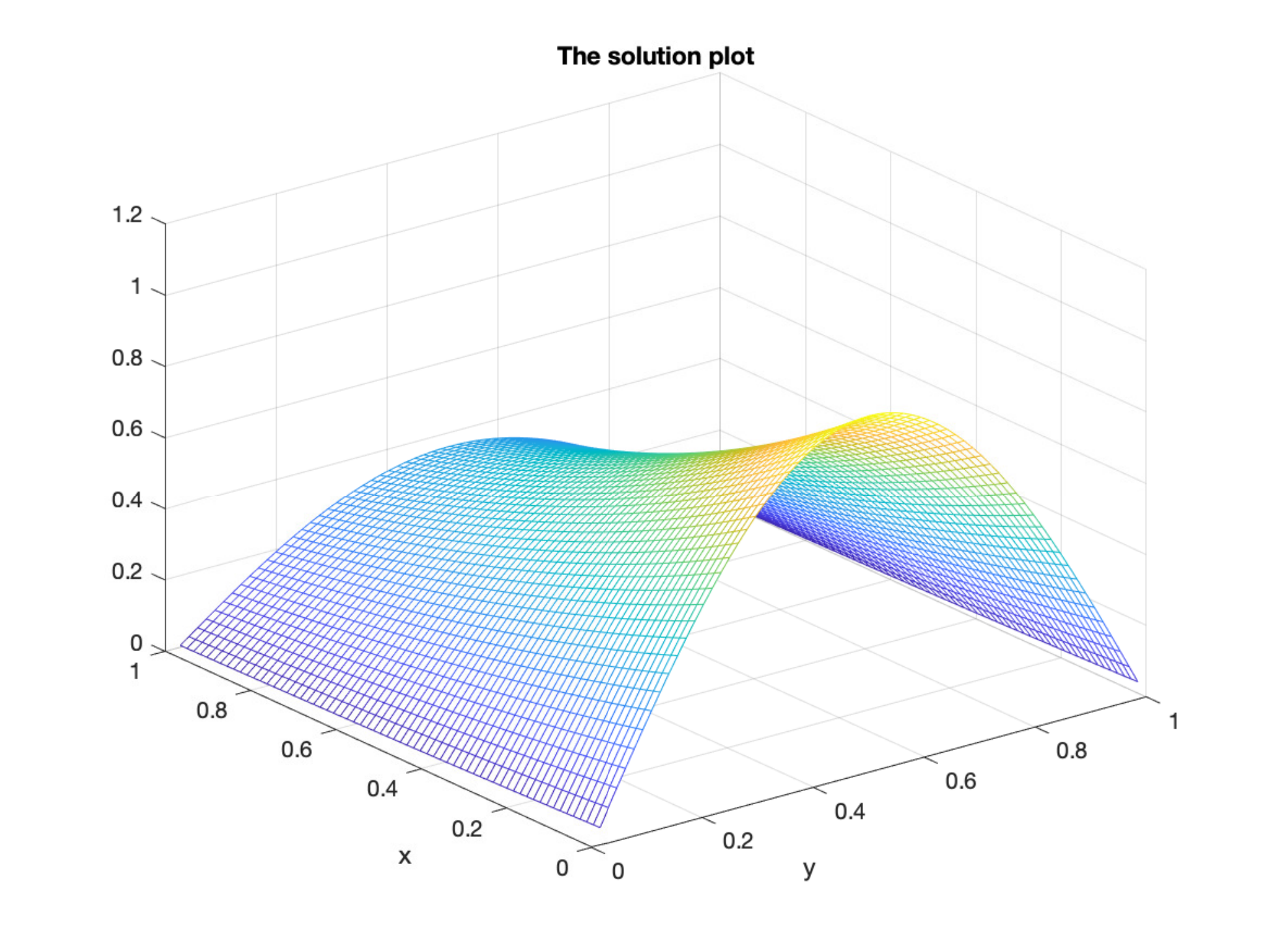}
 \end{minipage}
\begin{minipage}[t]{3.0in} (b)

\includegraphics[width=1.0\textwidth]{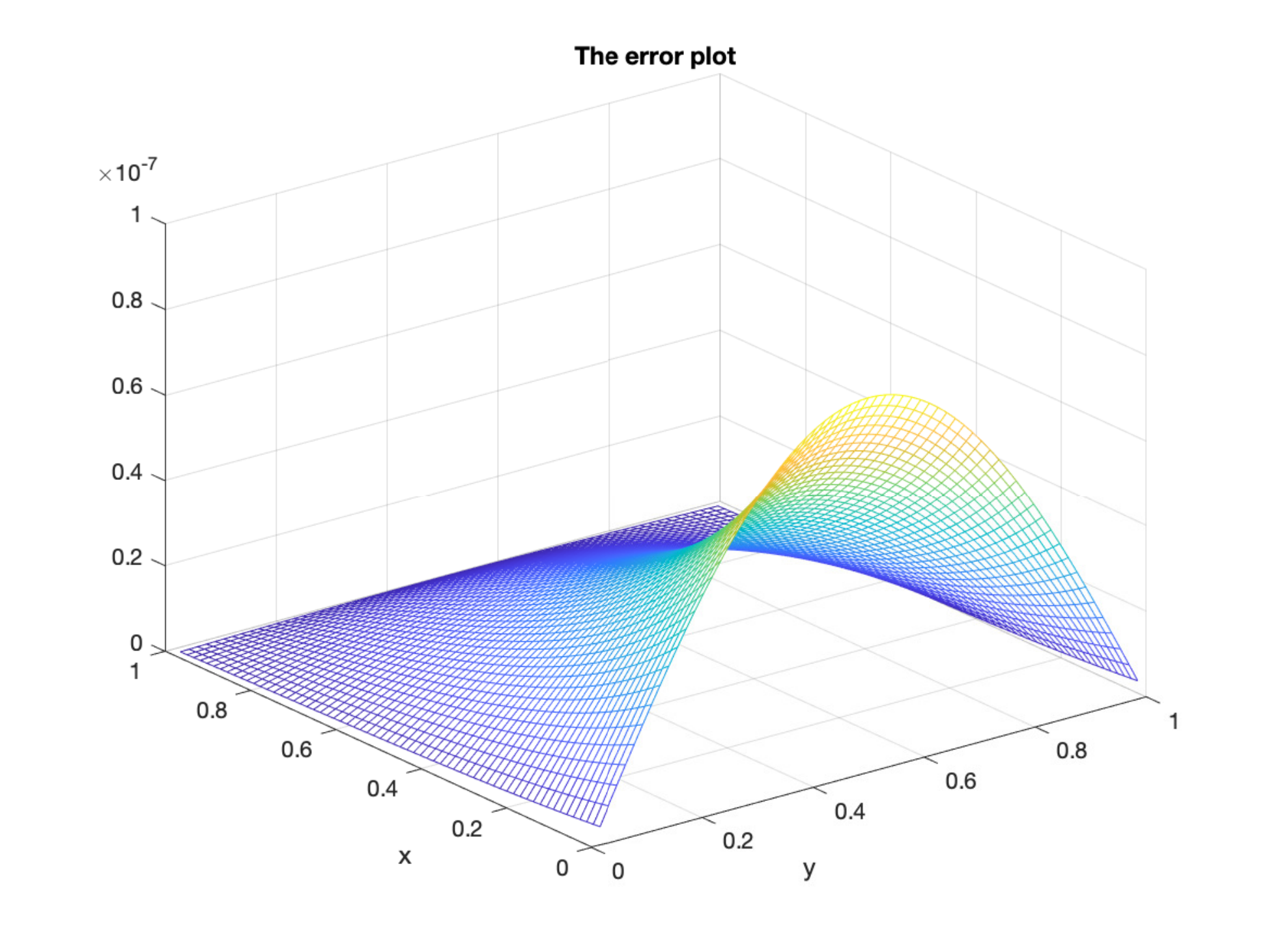}
 \end{minipage}

\caption{(a): A plot of the computed solution of Example~\ref{ex1a}  for the  Poisson  equation   and a Neumann BC  obtained from a $64$ by $64$ grid. (b): The error plot  where  $\|E\|_{\infty}=8.7602\times 10^{-8}.$}
\label{soln-err-plot-ex1}
\end{figure}

In Figure~\ref{soln-err-plot-ex1},  we show the solution and error  plots of  Example~\ref{ex1a}  for  the Poisson equation with  a Neumann BC obtained using a $64$ by $64$ grid. Both the solution and error are smooth and the error  is small,  $\|E\|_{\infty}=8.7602\times 10^{-8}.$

\begin{figure}[phbt]
\begin{minipage}[t]{3.0in} (a)
\includegraphics[width=1.0\textwidth]{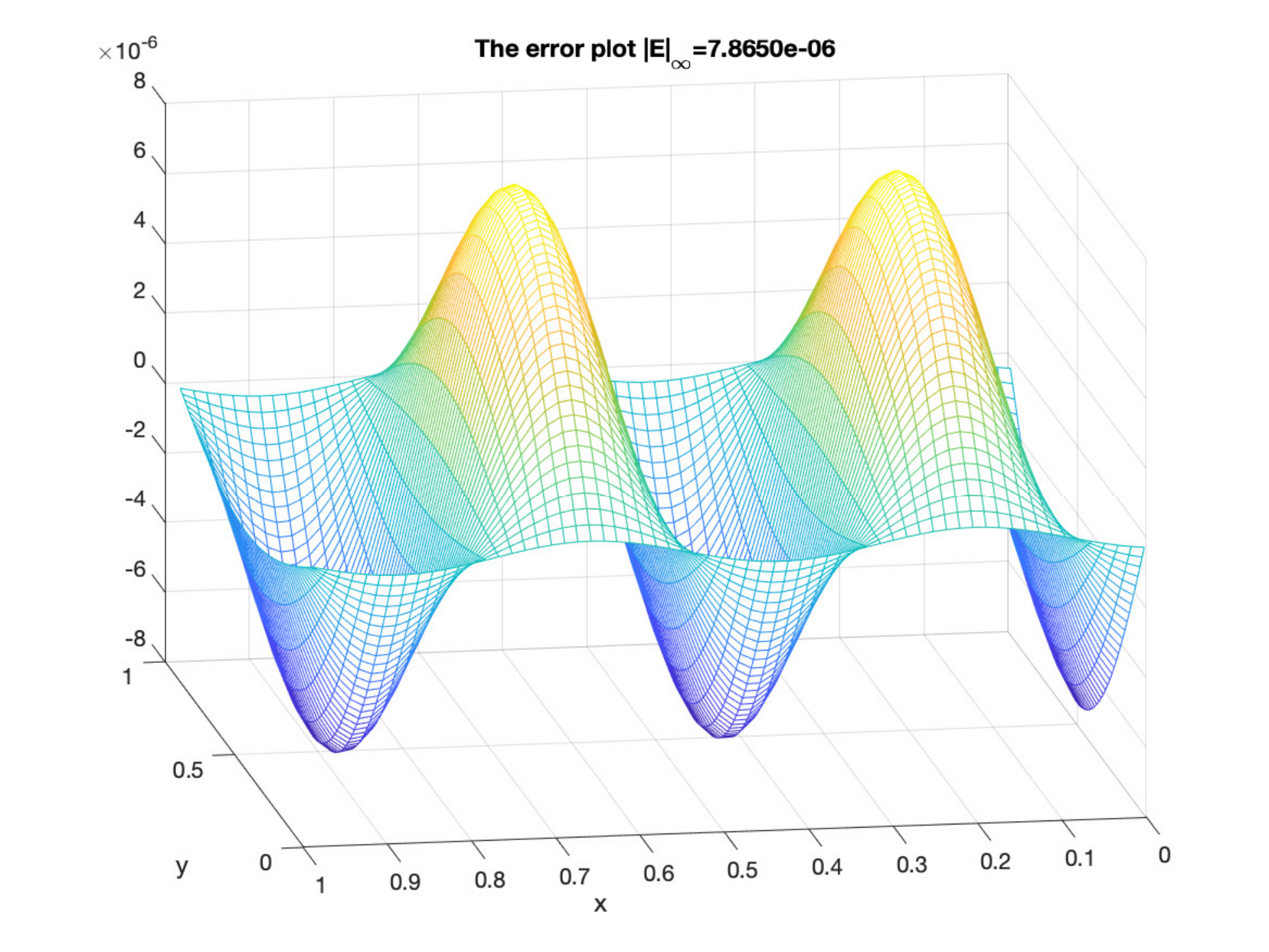}
 \end{minipage}
\begin{minipage}[t]{3.0in} (b)
\includegraphics[width=1.0\textwidth]{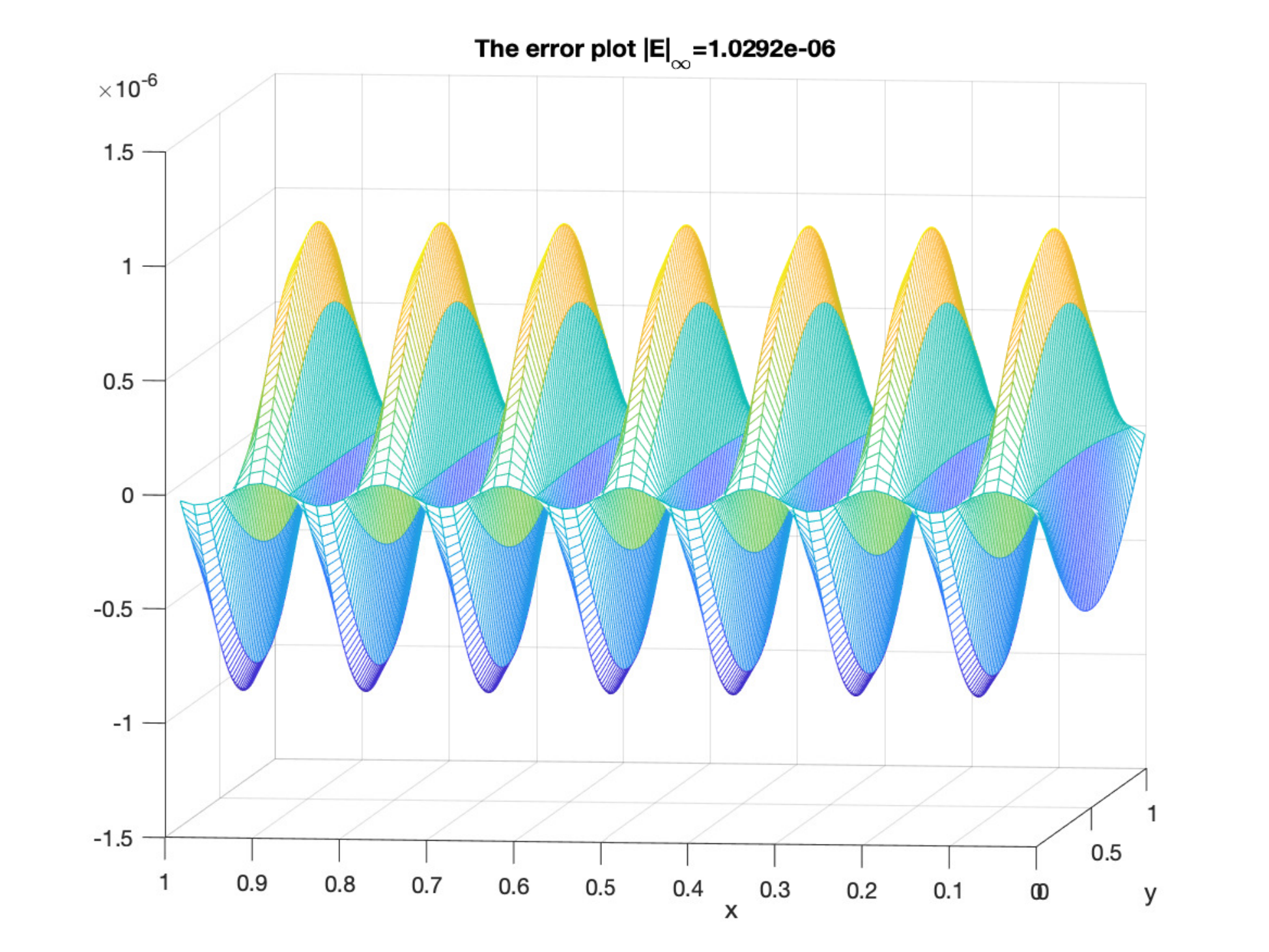}
 \end{minipage}

\caption{(a): Error plots of the computed solution using the compact fourth order schemes  for Example~\ref{ex1a} for the Helmholtz equation. (a): $K=200$ and the error has mild oscillation.   (b):  $K=2000$,  the error is more oscillatory. }
\label{err-plot-helm-ex1B}
\end{figure}

 In Figure~\ref{err-plot-helm-ex1B}~(a), we show an error  plot   from a $64$ by $64$ grid  of Example~\ref{ex1a} for the  Helmholtz  equation with $K=200$, $k_1=5$, $k_2=50$,  and a Neumann BC. The error now has mild oscillations.
In Figure~\ref{err-plot-helm-ex1B}~(b), we show  the error plot when the wave number is relatively large $K=2000$.   We see the error is  oscillatory even though the solution is smooth.

\begin{figure}[phbt]
\begin{minipage}[t]{3.0in} (a)

\includegraphics[width=1.0\textwidth]{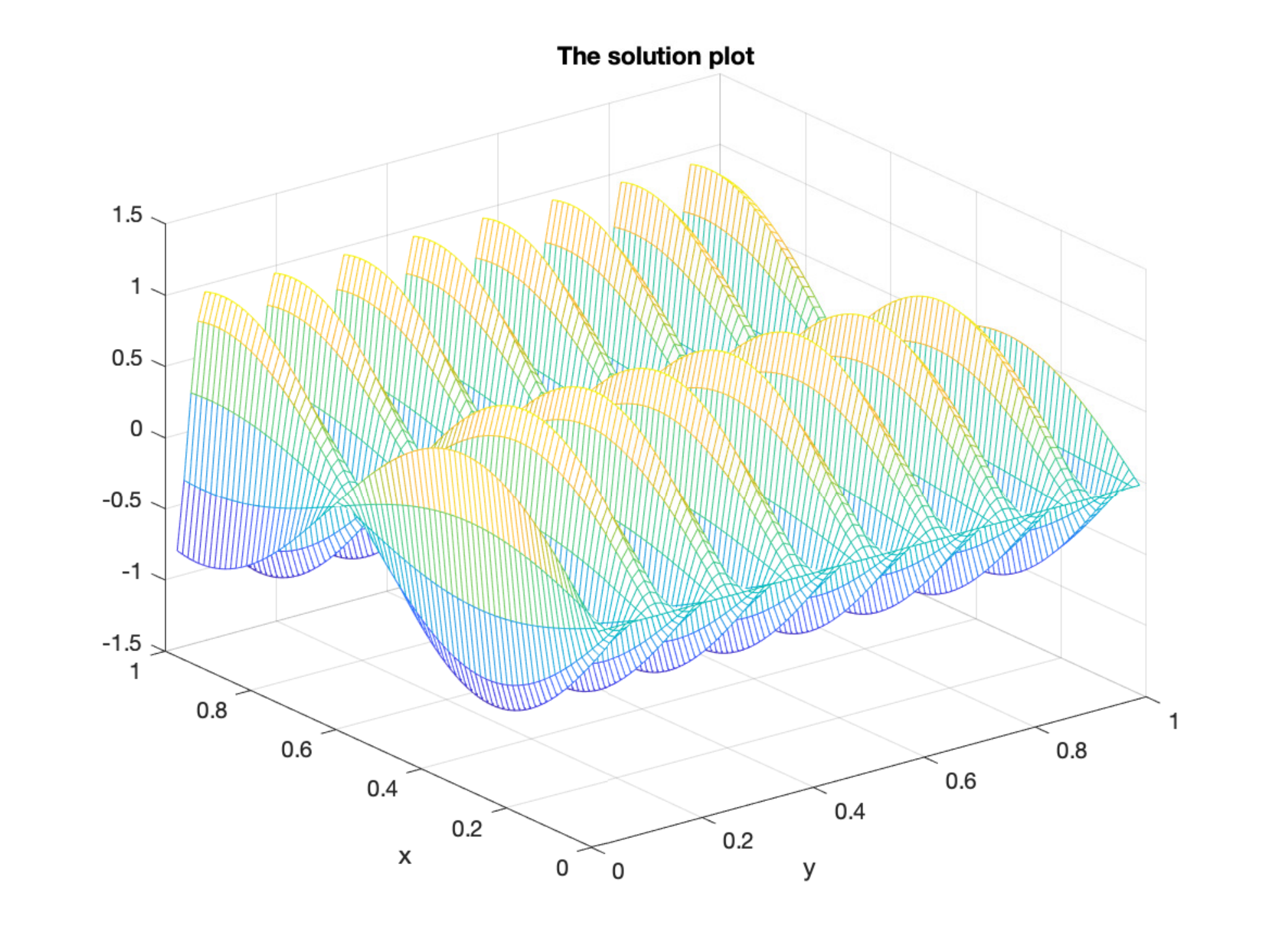}

 \end{minipage}
\begin{minipage}[t]{3.0in} (b)

\includegraphics[width=1.0\textwidth]{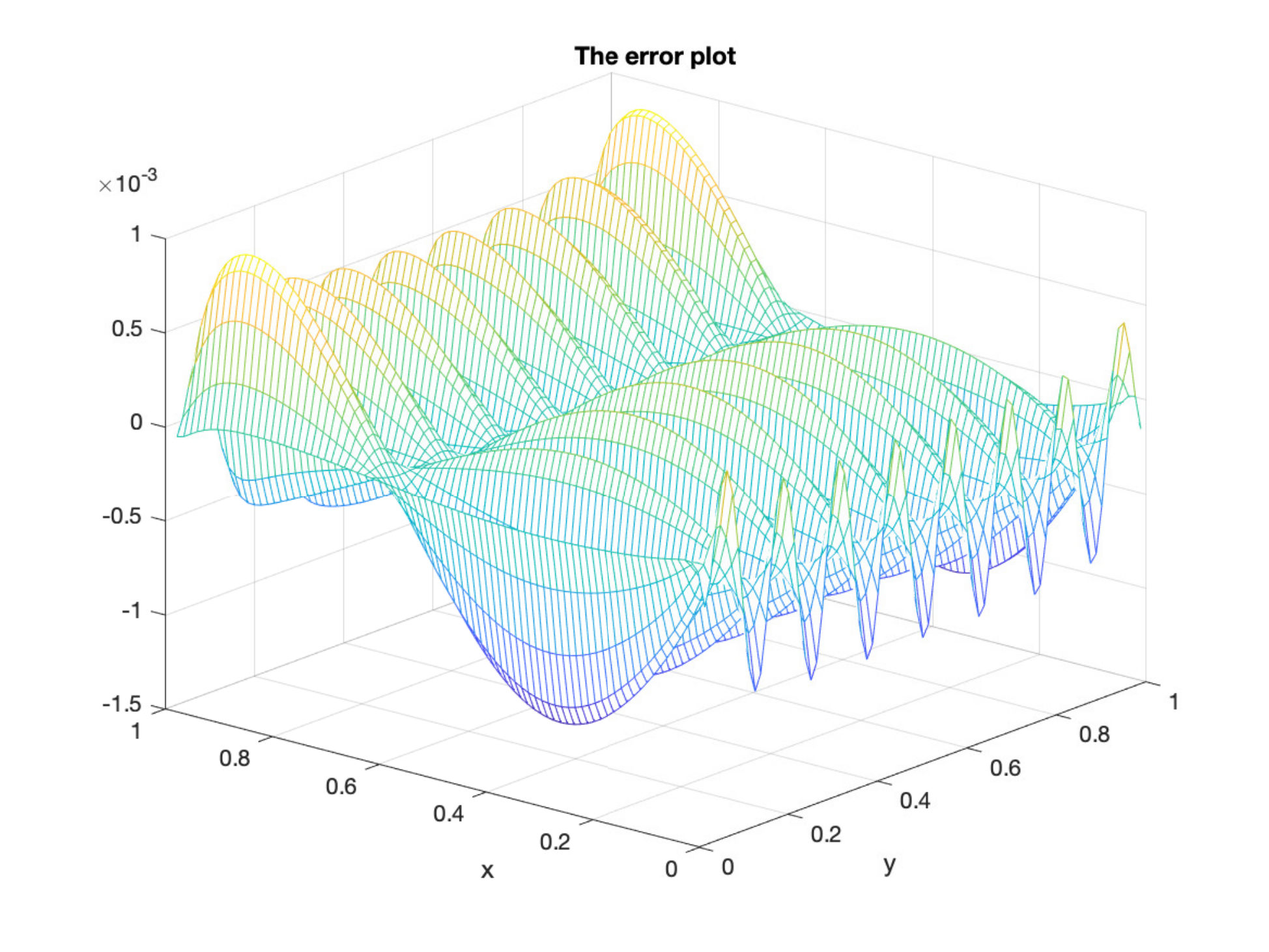}
 \end{minipage}
\caption{(a): The computed solution using the  fourth order compact schemes for Helmholtz equations  applied to Example~\ref{ex2} with
 $K=2000$, $k_1=5$, $k_2=50$.   (b): The error plot ($\|E\|_{\infty}=2.7619\times 10^{-3}$). Both  the solution and the error are oscillatory with some boundary effect $x=0$ where a Robin BC is specified.}
\label{err-plot-helm-ex2}
\end{figure}

 In Figure~\ref{err-plot-helm-ex2}~(a), we show a solution plot of Example~\ref{ex2}  for the  Helmholtz  equation with  $k_1=5$, $k_2=50$   from a $64$ by $64$ grid  and a Neumann boundary condition at $x=0$ and the Dirichlet  boundary condition elsewhere. The solution  is oscillatory.
In Figure~\ref{err-plot-helm-ex2}~(b) we show  the error plot when the wave number is  $K=2000$.   We see that the  error  is  also  oscillatory, and some boundary effect at $x = 0$.

\section{Super-third order compact schemes for flux BCs without $f$-extension} \label{sec:super-3rd}

If we do not use any extension of $f(x,y)$, then  we can seek a high order compact scheme of the following form
\eqml{3rd-hoc}
&  \dsp  \sum_{i_k=0}^{1} \sum_{j_k=-1}^{1} \!\!  \alpha_{i_k,j_k} U_{i+i_k,j+j_k} =  \sum_{i_k=0}^{1} \sum_{j_k=-1}^{1}  \!\!  \beta_{i_k,j_k} f(x_{i+i_k},y_{j+j_k})  +  \sum_{j_k=-1}^1 \!\! \gamma_{j_k}   g(y_{j+j_k})  , \\ 
  & \dsp  \sum_{i_k=0}^{1} \sum_{j_k=-1}^{1}  \!\!  \beta_{i_k,j_k} = 1.
\enml
Now the indexes for $f_{ij}$ and $U_{ij}$ range from $i=0,1$ and $j=-1,0,1$. The degree of freedom  is $15$ which is not  enough for a fourth-order scheme.  
Another consideration is the stability. We want the coefficient matrix is an M-matrix if $K\le 0$.

It is important to have both consistency and stability. Thus, we  give up two equations corresponding to $x^4$ and $y^4$ while keep other equations. In the modified system of equations, we will have HOC schemes that  are better than third  but not fully fourth order accurate. Therefore, we call such schemes  super third-accurate methods.

{For the stability concern, we use a maximum principle preserving scheme to enforce the sign property.
Let the system of linear equations for the coefficients be $A {\bf x } = {\bf b}$, where $A\in R^{14\times 15}$.
We impose the sign restrictions on the coefficients
$\{ \alpha_{i_k,j_k} \}$ in \eqn{3rd-hoc}
\eqml{sign}
 \alpha_{i_k,j_k}   \ge 0 \quad \mbox{if} \quad (i_k,j_k) \ne (0,0), \\ \eqsp
 \alpha_{i_k,j_k} <0  \quad \mbox{if} \quad (i_k,j_k) = (0,0),
\end{array}\end{equation}
along the equality constraints.

We form the following quadratic constrained optimization problem to determine
the coefficients of the finite difference scheme, see for example \cite{li:book},
\eqm \label{op1}
&& \dsp \min_{\bf  x } \left \{
 \half { \bf x}^T H {\bf x } - {\bf x  }^T {\bf w} \right \},\\ \eqsp
 && s.t.  \quad \left \{ \begin{array}{l}
    \dsp A {\bf x}  = {\bf b} \\ \eqsp
\dsp  \alpha_{i_k,j_k}  \ge  0, \quad \mbox{if} \quad (i_k,j_k) \ne (0,0), \\ \eqsp
\dsp  \alpha_{i_k,j_k}  < 0, \quad \mbox{if} \quad (i_k,j_k) = (0,0),  
\end{array} \right. \label{op2}
\enm
where ${\bf x}$ is the vector composed of
the coefficients of the finite difference equation,    the coefficients of the combination of $f_{ij}$, and  the coefficients of the combination of
$g(y_{j_k})$. In the implementation, we take $H=A^T A$,
 and ${\bf w} = A^T {\bf b}$, and  use the Matlab quadratic programming function `quadprog' to solve the optimization problem with the initial guess  ${\bf x}_0 = A^+ {\bf b}$, where $A^+$ is the pseudo-inverse of $A$. It is possible to have better $H$ and $\bf w$.

\begin{figure}[phbt]
 \centering
 \includegraphics[width=0.415\textwidth]{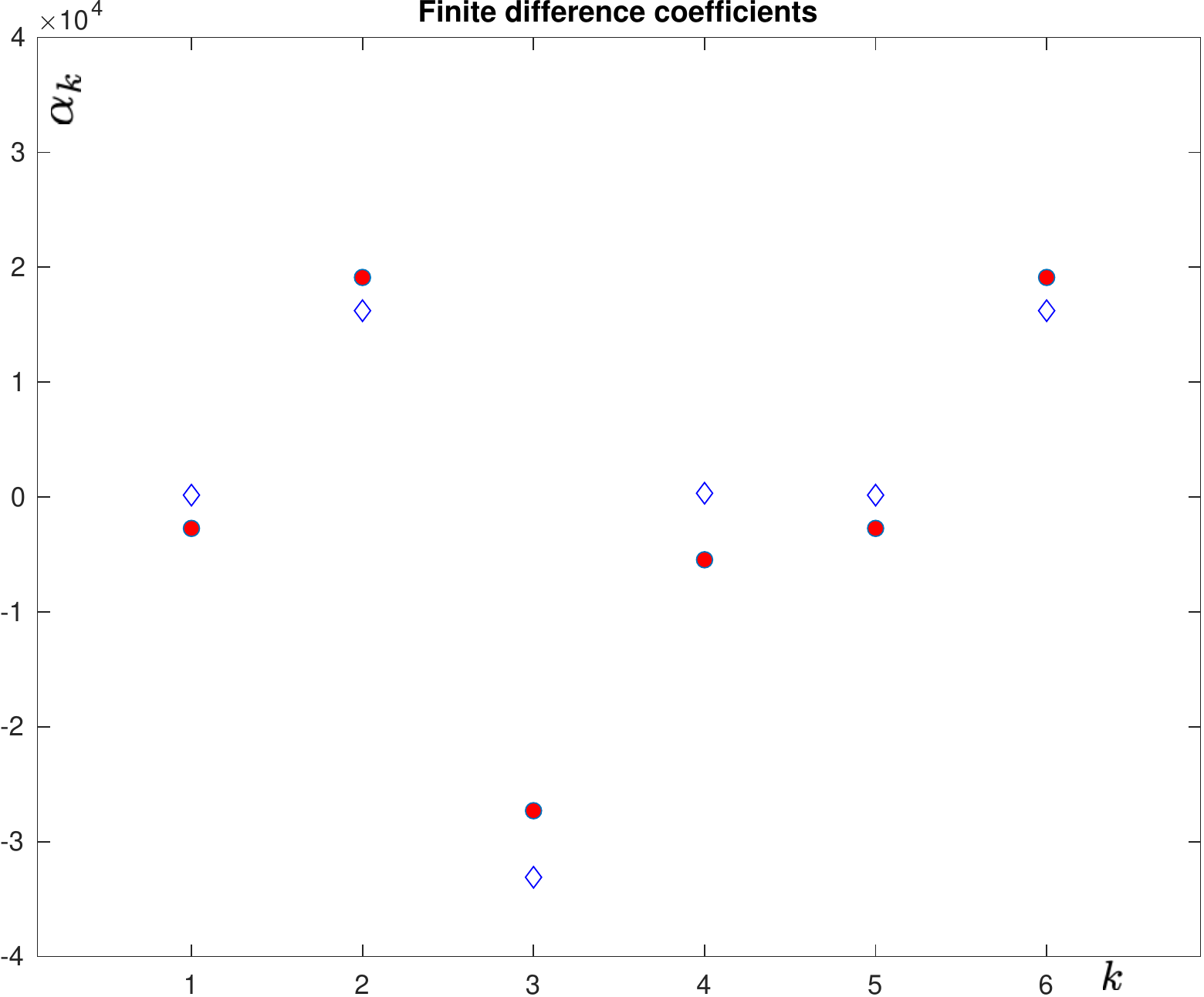}
\caption{The finite difference coefficients for $U_{0j}$ without (the filled red circles);  and with the optimization process (the blue diamonds).  The coefficient matrix of the FD equations using the second one is an M-matrix. The six coefficients $\alf_k$'s correspond to the grid points
$(0,y_j-h), (h, y_j-h), (0,y_j), (h,y_j), (0, y_j+h), (h,y_j+h)$,  where $(0,y_j)$ is the master grid point. }   \label{opt-coe}
\end{figure}

 In Figure~\ref{opt-coe} we show a computed set of the finite difference coefficients for $U_{0j}$ at a grid point, say $(x_0,y_j)$ where a  Neumann or Robin boundary condition is prescribed. The third coefficient is the one on the diagonal. The data marked with the filled red circles is obtained directly from ${\bf x}_0 = A^+ {\bf b}$ without optimization, where $A^+$ is the pseudo-inverse of the coefficient matrix $A$. We can see that the coefficients of the diagonals can be positive or negative. The data marked with  the blue diamonds  is obtained  from the optimization process.  We can see that the coefficients off the diagonals now are all non-negative.}

Once again, there are infinite number of solution sets to the linear system of equations. We list one particular ideal set below,
\begin{align} \label{plus3-coeff}
U_{i,j}:  \;   \frac{1}{\lambda h^2}\left[
               \begin{array}{cc}
                  4 - Kh^2 & 2 \\
                   2\sigma Kh^3 + (\lambda + 2)Kh^2 - 12\sigma h - 20 \null \quad & 8 \\
                  4 - Kh^2 & 2 \\
               \end{array}
             \right],
\end{align}
\begin{align} \label{plus3-coeffB}
 f_{i,j}: \; \frac{1}{\lambda}\left[
               \begin{array}{cc}
                  0 & 0 \\
                  4-Kh^2 \null \quad & 2 \\
                  0 & 0 \\
               \end{array}
             \right],   \quad
g_{j}: \;   \frac{1}{h}\left[
               \begin{array}{c}
                  0  \\
                  -2  \\
                  0  \\
               \end{array}
             \right],
\end{align}
where $\lambda = 6-K h^2$.  Note that, if $K\le0$, $\sigma\ge 0$,  then the coefficient matrix is an M-matrix and the HOC  scheme preserve the discrete maximum principle.


Similar to the error estimates in the previous section, we have the convergence theorem for the super-third compact  method for the particular set of coefficients list above as stated below.
\begin{theorem}
Let $U_{ij}$ be the finite difference solution obtained from the super-third order compact scheme with a set of coefficients given by (\ref{plus3-coeff})-\eqn{plus3-coeffB}. Assume that the solution to the Poisson/Helmholtz
equation with non-empty ${\cal{R}}_1$(Dirichlet) and  ${\cal{R}}_2$(Robin) in the boundary condition~\eqn{two-BC} is $u(x,y)$. Then, the algorithm is {\em exact}  if the solution is any fourth  order polynomials without $x^4$ and $y^4$ terms when $K=0$ and $\sigma=0$. For general solutions,  assuming  $u(x,y)\in C^6(\cal{R})$,  we have the following error estimate
\eqm
  \left \| u(x_i,y_j)- U_{ij} \right \|_{\infty}  \le C \, max\{ \|D^6 u\|_{\infty, {\cal{R}}},  \, \|D^5 u\|_{\infty, \partial {\cal{R}}}\}  \, h^{3+},
\enm
\end{theorem}
 The proof is similar to that of Theorem~\ref{thm-A}.
 From the PDE theory we know that, if $K\le 0$ and $\sigma\ge 0$, then the continuous problem is wellposed, while $K>0$ is not guaranteed.


\vthin

\begin{table}[htbp]
\caption{Grid refinement analysis  of the super-third order compact scheme for  Helmholtz  equations with a Neumann BC. (a):  Example~\ref{ex1a} with $K=200$. The average convergence order is $3.3118$.  (b): Example~\ref{ex2} with $K=2000$, $k_1=5$, $k_2=50$. The convergence order is clean fourth order. (c):   Example~\ref{ex2} with $K=2000$, $k_1=25$, $k_2=5$.  The average for the last column is $4.0393$.} \label{third_plus_Helm}

\vthin


\begin{center}
\begin{tabular}{|c| c c |}
  \hline
       $N$   & $\|E\|_{\infty} $  & $order$    \\       \hline	
   $8$  &  7.9660\,e-03    & \\
  $16$  &  5.6584\,e-04  &  $3.8154$\\
  $32$  &  6.5175\,e-05  & $3.1180 $ \\
  $64$  &  7.8650\,e-06  & $ 3.0508 $ \\
 $128$  & 1.2122\,e-06   & $2.6978$\\
 $256$  & 1.4943\,e-07  & $3.0201$  \\
 $512$  &  8.3095\,e-09  & $ 4.1686$ \\
      \hline
\end{tabular} \; $\null \quad $
\begin{tabular}{|c| c c | }
  \hline
       $N$   & $\|E\|_{\infty} $  & $order$    \\       \hline	
  $16$  &  1.0922\,e-00  &  $$\\
  $32$  &  4.7754\,e-02  & $  4.5138$ \\
  $64$  &  2.7660\,e-03  & $ 4.1114 $ \\
 $128$  & 1.7023\,e-04   & $4.0222$\\
 $256$  & 1.0526\,e-05  & $ 4.0155$  \\
 $512$  &  6.6300\,e-07  & $  3.9888 $ \\
      \hline
\end{tabular}\; \; $\null \quad $
\begin{tabular}{|c| c c | }
  \hline
       $N$   & $\|E\|_{\infty} $  & $order$    \\       \hline	
  $16$  &  	0.5620\,e-00  &  $$\\
  $32$  &  	1.1471\,e-02  & $ 5.6145 $ \\
  $64$  &  7.7896\,e-04  & $ 3.8803 $ \\
 $128$  & 4.1507\,e-05   & $4.2301$\\
 $256$  & 	2.7972\,e-06  & $ 3.8913$  \\
 $512$  &  1.5868e-07\,e-07  & $   4.1398$ \\
      \hline
\end{tabular}
\end{center}
\end{table}

We have repeated numerical tests for the same examples in the previous section  without extension $f(x,y)$ using the super-third compact method with a  flux type BC.
In Table~\ref{third_plus_Helm}, we show grid refinement analysis  results  for the  Helmholtz  equation.
In Table~\ref{third_plus_Helm}~(a), we list the results  for  Example~\ref{ex1a} with $K=200$. The order of convergence fluctuates between three and four and the average  is $3.3118$. In Table~\ref{third_plus_Helm}~(b), the test is for  Example~\ref{ex2} with $K=2000$, $k_1=5$, $k_2=50$.  The results show a  fourth order convergence. One of explanations is that the error is dominated in the $y$ direction ($k_2=50$) rather than from the Neumann BC at $x=0$. In Table~\ref{third_plus_Helm}~(c), the test is for  Example~\ref{ex2} with $K=2000$, $k_1=25$, $k_2=5$.   The order of convergence fluctuate between three and four and the average  is $4.0393$, which shows some natures of the super-third compact method.


\section{HOC schemes for diffusion and advection equations} \label{sec:diff-advec}

In this section, we show that the same  idea  can be applied to diffusion and advection equations with constant coefficients,
\eqm
    \Delta u + a u_x + b u_y + Ku = f.
\enm

\ignore{
\textcolor{blue}{
The fourth order compact scheme of the above convection-diffusion equation with $K=0$ from \cite{MR1440337} for interior grid points is given by
\begin{equation}
  \sum_{i_k=-1}^{1} \sum_{j_k=-1}^{1} \!\!  \alpha_{i_k,j_k} U_{i+i_k,j+j_k}  =  \frac{1}{12}\left(f_{i-1,j} + f_{i+1,j}+f_{i,j-1} + f_{i,j+1} + 8 f_{i,j}
  +\gamma(f_{i+1,j}-f_{i-1,j})+\delta(f_{i,j+1}-f_{i,j-1})\right),
\end{equation}
where $\gamma=ah/2$ and $\delta =bh/2$ are the cell Reynolds numbers, and $\alpha_{i_k,j_k}$ is depicted by the nine-point stencil shown in Figure \ref{stencil}.
If $a=b=0$, the above scheme reduces to the classical compact difference scheme (\ref{compact2d}) for 2D Poisson equation.}

\begin{figure}[phbt]
 \centering
  \includegraphics[width=0.6\textwidth]{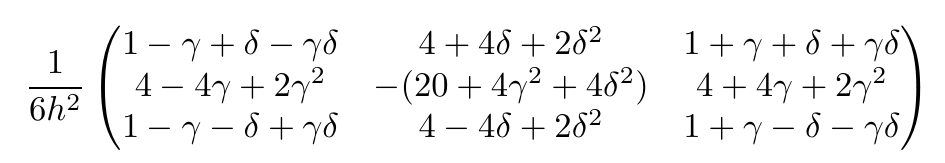}
\caption{\textcolor{blue}{Compact nine-point stencil.}} \label{stencil}
\end{figure} }

Using the same approach, the finite difference equation for an {\em interior grid point} \ can be written as
\eqm \label{hoc-diff-advec}
&&  \dsp  \sum_{i_k=-1}^{1} \sum_{j_k=-1}^{1} \!\!  \alpha_{i_k,j_k} U_{i+i_k,j+j_k}  =  \sum_{i_k=-1}^{1} \sum_{j_k=-1}^{1}  \!\! \beta_{i_k,j_k} f(x_{i+i_k},y_{j+j_k})  , \\ 
  && \dsp  \sum_{i_k,j_k=-1}^1 \!\!\! \beta_{i_k,j_k} = 1. \label{hoc-diff-advecB}
\enm
The coefficients are determined in a similar way as discussed in previous sections. We list additional PDE relations and the system of linear equations for the coefficients in Appendix. There are infinite number of solutions of the coefficients.  With the Matlab/Maple symbolic package, we have obtained  a set of solution  of coefficients for the diffusion and advection PDEs that has the following form
\eqml{coe-split}
 & \dsp  {\bfalf} = \frac{A_{diff} }{h^2} +  \frac{B_{advec}(a,\,b) }{h} + C_{other}(K), \\ \eqsp
 & \dsp A_{diff} \sim O(1), \qquad B_{advec}(a,\,b) \sim O(a,b), \qquad C_{other} \sim O(K),
\enml
where $\frac{A_{diff} }{h^2}$ is the leading terms in the FD coefficients which has the symmetry in terms of the four corner  and  west-east-north-south grid points in reference to the master grid point.
When $a=0$ and $b=0$, $\frac{A_{diff} }{h^2}$  is the same as that listed  in (\ref{Coeff_H_A})-\eqn{Coeff_H_B}.

The FD coefficients corresponding to $U_{ij}$ of the HOC scheme  for diffusion and advection equations computed with the Maple are shown in Figure~\ref{FD_Coef_MapleB} in Appendix.
The solution set is obtained using the  pseudo-inverse of the coefficient matrix, also called the SVD solution. The computation can be done almost instantly with  numerical solutions but took quite a while to return the symbolic (analytic) solution. As expected, the SVD solution is often the best compared with other set of solutions even if the convergence order is the same. We have seen that the error constant  can be several order magnitude smaller using the SVD solution than that using other least squares solutions.
In terms of the convergence proof, we can see that  all the terms that contain $a$ and $b$ are at most $O(1/h)$. The coefficient matrix is an M-matrix if $h$ is small enough, which leads to asymptotic fourth order convergence. In Figure~\ref{FD_Coef_MapleC}, we also list the combination coefficients of $f_{ij}$'s whose sum equals  one.

If a Robin or Neumann boundary condition, say, is defined at $x=x_l$, then we find {\em another set of FD coefficients} of the form as in \eqn{4th-hoc}  with different PDE, high order PDE relations, and resulting linear system of equations for the coefficients assuming that we have an extension of $f(x,y)$. We also list those relations in Appendix. The FD coefficients for $U_{ij}$  also have the form,
\eqm
  {\bfalf}_{Robin} = \frac{\bar{A}_{diff} }{h^2} +  \frac{\bar{B}_{advec}(a,\,b,\sigma) }{h} + \bar{C}_{other}(K).
\enm
The coefficient matrix of the FD  is an M-matrix unless extreme situations when $K$, $|a|$, $|b|$ is so  large that the optimization process fails to return a feasible solution. The local truncation error is of $O(h^3)$ which would not affect the global fourth order convergence.

\subsection{Convergence analysis of the HOC method for diffusion and advection equations}

The discussion of the  convergence for the developed HOC scheme is challenging since it depends on the advection  coefficients and boundary conditions. Using the computed sets of coefficients, and the fact that the local truncation errors at boundary grid points can be one order lower than that of interior grid points without affecting the global accuracy of the computed solution, we can obtain asymptotic fourth order convergence theorem.

\begin{theorem}
Let $U_{ij}$ be the finite difference solution obtained from the derived  HOC scheme for the diffusion and advection equation with non-empty ${\cal{R}}_1$(Dirichlet) and  ${\cal{R}}_2$(Robin) in the boundary condition~\eqn{two-BC}. Assume that the solution  $u(x,y)\in C^6(\cal{R})$ , if $h$ is  small enough,  $K\le 0$, and $\sigma\ge 0$, then the following error estimates hold
\eqm
  \left \|  T_h  \right \|_{\infty}  \le \left \{ \begin{array}{ll}
  \dsp   C_1 \, \|D^6 u\|_{\infty, {\cal{R}}}  h^4 &  ~ \mbox{interior grid points, } \\ 
  \dsp  C_2\,  \|D^5 u\|_{\infty, \partial {\cal{R}}} h^3  &    \mbox{flux boundary grid points,}
 \end{array} \right.   \qquad    \left \| u(x_i,y_j)- U_{ij} \right \|_{\infty}  \le {\bar C}   h^4,
\enm
for Dirichlet, or Dirichlet with part of Robin, or Neumann boundary conditions.
\end{theorem}

{\bf Proof:} At an interior point,   
we carry out the Taylor expansion at a grid $(x_i,y_j)$ of the local truncation error for all terms ($u(x_i,y_j)$ and $f(x_i,y_j)$ involved.  Thus, we should have
\eqm
  T_{ij}^h = T_{diff} h^3 + T_{advec}(a,\,b) h^4 + C_{other}(K) h^4 + O(h^5) ,
\enm
according to \eqn{coe-split}. Furthermore, since $A_{diff}$ is symmetric, or centered discretization in reference $(x_i,y_j)$, the terms in the expansion involving odd partial derivatives are canceled out,  which leads to  $\left \|  T_h  \right \|_{\infty} \sim O(h^4)$ if $h$ is small enough.
Thus, for a Dirichlet boundary condition, the HOC scheme is asymptotically fourth order convergent.

If part of flux boundary condition is given,  then there  are no cancellations from $1/h^2$ terms, thus, we have $\left \|  T_h  \right \|_{\infty} \sim O(h^3)$. 
Since the coefficient matrix is an M-matrix which is  true if $K\le 0$,    $h\le  C/\max\left \{|a|, |b|\right \}$, and $\sigma \ge 0$,
 we  apply the convergence theorem in \cite{morton-mayers} to get the asymptotically fourth order convergence.  \hfill $\square$

\subsection{Numerical examples  of the HOC method for diffusion and advection equations}

We carried out  numerical experiments for Example ~\ref{ex1a}  and Example~\ref{ex2} using the resultant source terms and the boundary conditions.
In Table~\ref{D_BC_diff_advec}, we show numerical experiments results and the grid refinement analysis.  The left and middle tables  are the results  of  the method applied to Example~\ref{ex1a} while the right table lists results  of  the method applied  to Example~\ref{ex2}.  In the left table, the parameters are $K=20$, $a=1$, $b=2$ in which the convection is not very large. The results show clearly fourth order convergence starting from a  rather coarse grid
$N=16$. In the middle of table, the  parameters are $K=20$, $a=100$, $b=5$ in which the convection is  relatively strong. The fourth  order convergence was affected at the  coarse grids level until $N\ge 128$.  In the right table,  the  parameters are $K=20$, $a=1$, $b=100$, $k_1=5$, $k_2=10$.  With modest $k_1$ and $k_2$, Example~\ref{ex2} is tougher to compute due to the oscillations and  larger magnitudes of the partial derivatives and the source term. Nevertheless,  when the grid is fine enough, we   see clearly fourth order convergence.
\begin{table}[htbp]
\caption{Grid refinement analysis of the fourth order compact scheme for the diffusion and advection equation with a Dirichlet BC. (a): Results for
Example~\ref{ex1a} with $K=20$, $a=1$, $b=2$. (b): Results for   Example~\ref{ex1a} with large convection $a=100$, $b=5$.  (c):   Results for  Example~\ref{ex2} with $K=20$, $a=1$, $b=100$, $k_1=5$, $k_2=10$.  } \label{D_BC_diff_advec}

\vthin

\vthin

\begin{center}
\begin{tabular}{|c| c c |}
  \hline
       $N$   & $\|E\|_{\infty} $  & $order$    \\       \hline	
 16 & 1.0772e-04 & \\
32 &	6.7160e-06 &  4.0035     \\
64	& 4.1993e-07 & 3.9994  \\
128	& 2.6249e-08 & 3.9998  \\
256	& 1.6666e-09 & 3.9773  \\
512	& 7.3411e-11 &4.5048 \\
      \hline
\end{tabular} \; $\null \quad $
\begin{tabular}{|c| c c | }
  \hline
       $N$   & $\|E\|_{\infty} $  & $order$    \\       \hline	
  16  &  7.2039e-05  &  $$\\
  32	& 8.3145e-06 & 3.1151  \\
64	&7.1788e-07 & 3.5338 \\
128	& 5.1167e-08 & 3.8105  \\
256	& 3.3067e-09 & 3.9517 \\
512	&1.6823e-10 & 4.2969    \\
      \hline
\end{tabular}\;  $\null \quad $
\begin{tabular}{|c| c c | }
  \hline
       $N$   & $\|E\|_{\infty} $  & $order$    \\       \hline	
  16	&1.0018e-02 &  \\
32&	1.5245e-03& 2.7162    \\
64	&1.4749e-04&3.3696 \\
128&	1.0652e-05&  3.7914 \\
256	&6.9281e-07& 3.9425 \\
512&	4.3801e-08 & 3.9834 \\
      \hline
\end{tabular}
\end{center}
\end{table}
\begin{figure}[phbt]
\begin{minipage}[t]{3.0in} (a)

\includegraphics[width=1.0\textwidth]{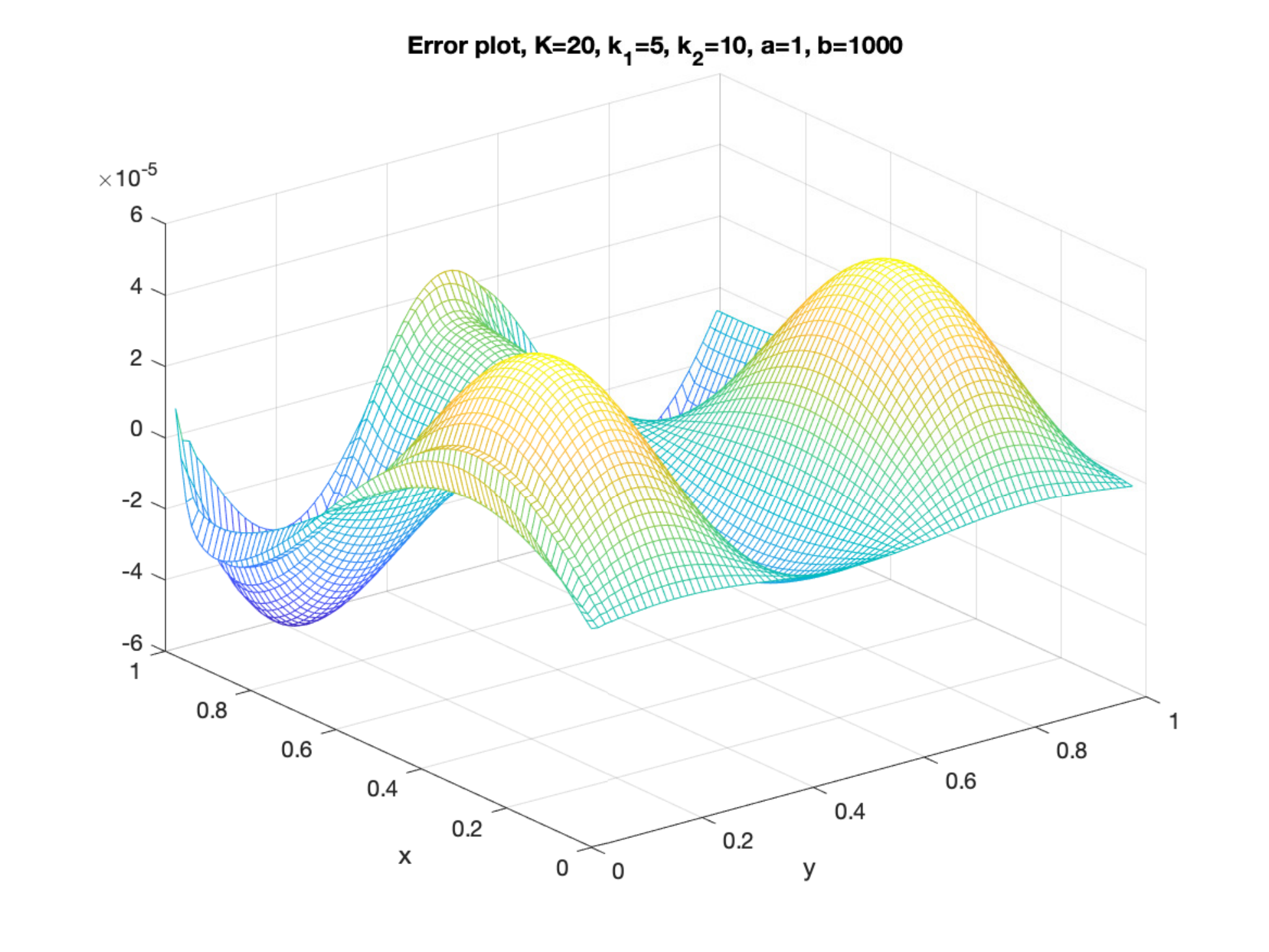}
 \end{minipage}
\begin{minipage}[t]{3.0in} (b)

\includegraphics[width=1.0\textwidth]{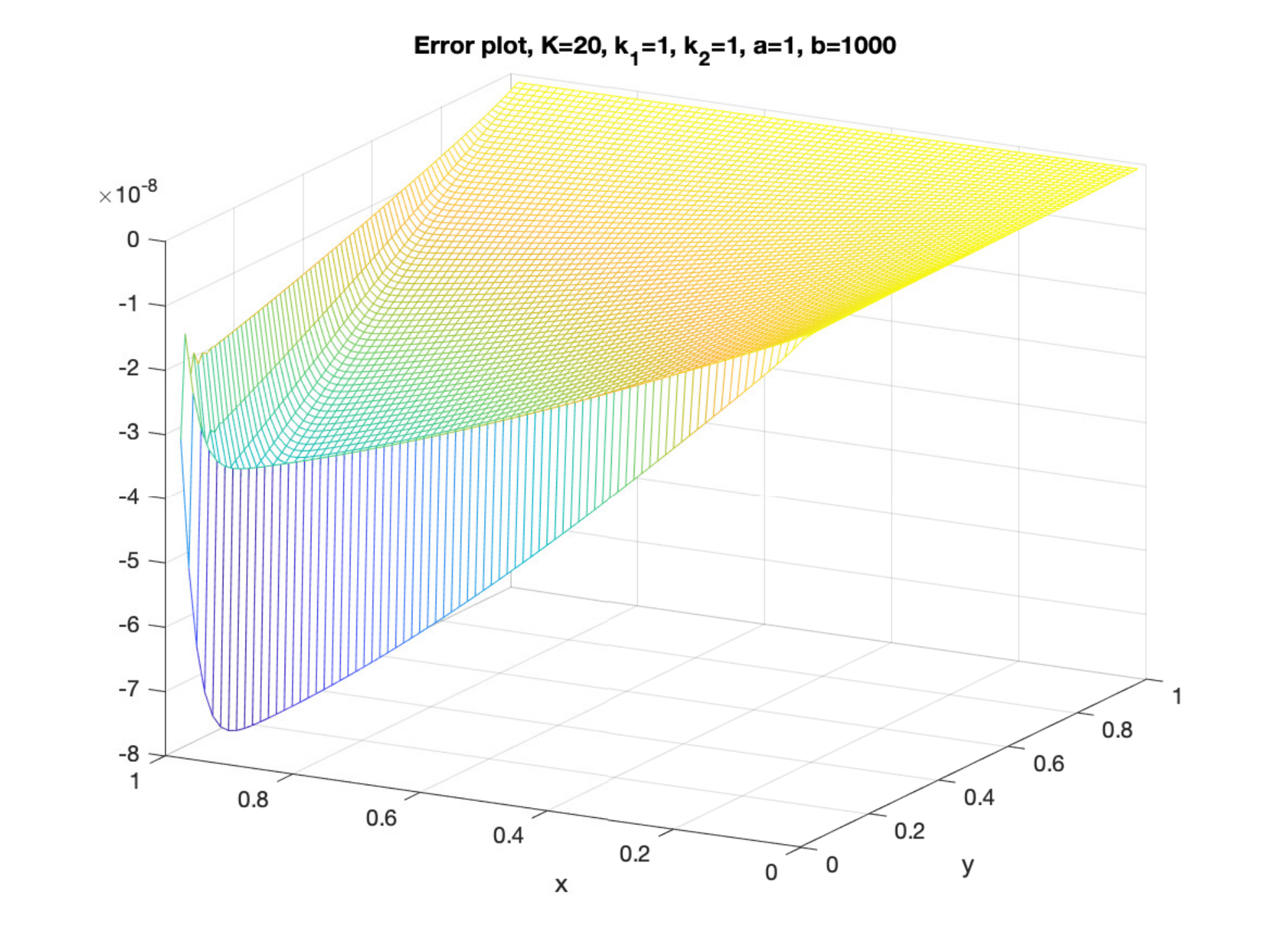}
 \end{minipage}
\caption{(a): Error plots of the computed solution using the  HOC schemes for  the diffusion and advection equation for Example~\ref{ex2}. (a): $K=20$, $a=1$, $b=10$, $k_1=5$, $k_2=1000$. The oscillation is dominant compared with the advection effect.   (b): $K=20$, $a=1$, $b=1000$, $k_1=1$, $k_2=1$.  We see the boundary layer effect in the error plot although the error ($10^{-8}$) is small.}
\label{err-plot-diffu-ex2}
\end{figure}

\ignore{
In Table~\ref{tab:diffu-ux}, we show  a grid refinement analysis  of the  HOC scheme with a Dirichlet boundary condition and larger oscillations.
  Due to the nature of the solution, the strong advection, large magnitude of the partial derivatives, which together leads to larger error constant. The mesh needs to be fine enough to observe  the asymptotic  fourth order convergence.  The  second-third  columns  in the table are the results for  Example~\ref{ex1a} with  $a = 100, b_1 = -10$, $K=20$. In Example~\ref{ex1a}, there is no oscillations in the solution and the errors are small even with coarse mesh. The  fourth-fifth columns  are the results for Example~\ref{ex2} with  $k_1=3$, $k_2=25$, and $a = 100, b_1 = -5$, $K=20$. In both cases, we see clearly asymptotic fourth order convergence. }

In Figure~\ref{err-plot-diffu-ex2}, we show two error plots with different parameters and  a strong  advection in  the $y$-direction for  Example~\ref{ex2} using a $80$ by $80$ grid. In both cases, the errors in the infinity norm are small, $10^{-5}$ and  $10^{-8}$. In Figure~\ref{err-plot-diffu-ex2}~(a), the oscillation is dominant  compared with the advection effect.  In Figure~\ref{err-plot-diffu-ex2}~(b), the solution is less oscillatory, we see the boundary layer effect clearly from the strong advection.


\ignore{
\begin{table}[hpbt]
\caption{Grid refinement analysis  of the HOC scheme for the diffusion and advection equation with a Dirichlet BC.  The  second-third  columns  are the results for  Example~\ref{ex1a} with  $a = 100, b_1 = -10$, $K=20$.   The  fourth-fifth columns  are the results for Example~\ref{ex2} with  $k_1=3$, $k_2=25$, and $a = 100, b_1 = -5$, $K=20$.  In both cases, we see clearly asymptotic fourth order convergence.} \label{tab:diffu-ux}
\begin{center}
\begin{tabular}{|c| c c | |c|c|}
  \hline
       $N$   & $\|E\|_{\infty} $  & $order$  & $\|E\|_{\infty} $  & $order$   \\       \hline	
  16	& 7.4325e-05 &  & 2.1024e-02 & \\
32	& 9.4213e-06 & 2.9798   & 1.6903e-03 & 3.6367 \\
64	& 8.9526e-07 & 3.3955 & 2.0523e-04 & 3.0420 \\
128	& 6.3879e-08 & 3.8089 & 1.5515e-05 & 3.7255 \\
256	& 4.1498e-09 & 3.9442  & 1.0152e-06 &  3.9338 \\
512	& 2.3371e-10 & 4.1503  & 6.5675e-08 & 3.9503   \\
      \hline
\end{tabular}
\end{center}
\end{table}
 }

In  Table~\ref{tab:diffu-robin},  we present numerical example for  Example~\ref{ex1a} with a Robin boundary condition. The  second-third  columns  are the results for  Example~\ref{ex1a} with  $K=50$,  and $a = 1, b_1 = -5$.
The errors are small even with coarse grids  since the advection is relatively small and we observe clean fourth order convergence.
 The  fourth-fifth columns  are the results  for  the same example  with  $K=50$, and $a = 5, b_1 = -100$, a relatively strong advection.
 We observe asymptotic fourth order convergence.
\begin{table}[hpt]
\caption{Grid refinement analysis  of the HOC scheme for a diffusion and advection equation with a Robin BC for Example~\ref{ex1a}.
 The parameters are $K=50$,  and $a = 1, b_1 = -5$ for the results in the second-third  columns, and  $K=50$, and $a = 5, b_1 = -100$ for the fourth-fifth columns. }\label{tab:diffu-robin}
 \begin{center}
 \begin{tabular}{|c| c c | |c | c|}
  \hline
       $N$   & $\|E\|_{\infty} $  & $order$    & $\|E\|_{\infty} $  & $order$  \\       \hline	
       16	& 4.5242e-04 &  & 7.2142e-05 & \\
32 &	2.8238e-05 & 4.0020 & 8.4424e-06 &  3.0951 \\
64	&1.7639e-06& 4.0008 & 7.7686e-07 & 3.4419 \\
128	&1.1020e-07 & 4.0006  & 5.5333e-08 & 3.8114 \\
256	& 6.9949e-09 & 3.9777 & 3.5909e-09 & 3.9457 \\
512	& 3.2561e-10 & 4.4251 & 1.8999e-10 & 4.2404 \\
      \hline
\end{tabular}
\end{center}
\end{table}

\begin{table}[hptb]
\caption{Grid refinement analysis  of the HOC scheme for a diffusion and advection equation with a Robin BC for Example~\ref{ex2}.
 The parameters are $K=50$, $k_1=5$, $k_2=50$,  and $a = 1, b_1 = -5$ for the results in the second-third  columns, and  $K=50$, $k_1=5$, $k_2=50$, and $a = 5, b_1 = -100$ for the fourth-seventh columns.   $\|T_h\|_{\infty}$  measures the local truncation errors.  The last two columns are the results with switched advection coefficients, $a=-100$ and $b_1=5$}.
\label{tab:diffu-robinB}
 \begin{center}
  \begin{tabular}{|c| c c | |c c || c c|}
  \hline
       $N$   & $\|E\|_{\infty} $  & $order$    & $\|E\|_{\infty} $  & $order$  & $\|T_h\|_{\infty}$ & $order$  \\       \hline	
       16	&5.9061e-04   &  &  5.5602e+03  & & 8.7058e00 & \\
32 &	3.7512e-05  & 3.9768 & 1.1095e-01  &  15.6129  & 1.2100e00 & 2.8470 \\
64	&2.3546e-06 &  3.9938 & 9.2102e-03  &  3.5905  & 1.2352e-01 & 3.2922 \\
128	&1.4734e-07 & 3.9983  &  6.3934e-04 & 3.8486 & 9.0471e-03 & 3.7711 \\
256	&  9.2020e-09 & 4.0011 & 4.1108e-05 & 3.9591 & 5.9032e-04 & 3.9379 \\
512	& 4.4578e-10 & 4.3675 &  2.5880e-06  & 3.9895 & 3.6881e-05 & 4.0005\\
      \hline
\end{tabular}
$\,$
  \begin{tabular}{||c c ||}
  \hline
       $\|E\|_{\infty} $  & $order$      \\       \hline
    1.4970e00 & \\
    6.1313e-02 &  4.6097 \\
    3.4927e-03 & 4.1338 \\
	2.1341e-04 &  4.0326 \\
	1.3891e-05 &  3.9414 \\
   8.7529e-07 &  3.9882 \\
      \hline
\end{tabular}       
\end{center}
\end{table}

In Table~\ref{tab:diffu-robinB}, we show the numerical result for a Robin boundary condition with $\sigma=20$ for  Example~\ref{ex2} (solution is oscillatory). The  second-third  columns  are the results with  $K=50$, $k_1=5$, $k_2=50$,
and $a = 1, b_1 = -5$. A clean fourth order convergence can be seen since the advection is relatively small.
 The  fourth-fifth columns  are the results  for  the same example  with  $K=50$, and $a = 5, b_1 = -100$, a relatively strong advection.
 We observe asymptotic fourth order convergence. The sixth-seventh columns  are the results  of the local truncation errors which also has an asymptotic fourth order.   The last two columns are the results with switched advection coefficients, $a=-100$ and $b_1=5$. We observe the similar convergence behavior. In both cases, the errors due to the convections are dominated than the errors from the flux  BC. 


 We also implemented the fourth order compact scheme from \cite{MR1440337} for interior grid points, and our  FD scheme  at   boundary grid points. The results are similar.
The fourth order compact scheme of the  convection-diffusion equation with $K=0$ from \cite{MR1440337} for interior grid points is given by
\eqml{gupta}
 \dsp  \sum_{i_k=-1}^{1} \sum_{j_k=-1}^{1} \!\!  \alpha_{i_k,j_k} U_{i+i_k,j+j_k}  &= & \dsp \frac{1}{12}\left( \frac{\null}{\null} f_{i-1,j} + f_{i+1,j}+f_{i,j-1} + f_{i,j+1} + 8 f_{i,j} \right. \\
  && \dsp \null \left. +\gamma \, (f_{i+1,j}-f_{i-1,j})+\delta(f_{i,j+1}-f_{i,j-1}) \frac{\null}{\null} \right),
\enml
where $\gamma=ah/2$ and $\delta =bh/2$, and $\alpha_{i_k,j_k}$  are the coefficients of the nine-point stencil shown in Figure \ref{stencil}. 

\begin{figure}[phbt]
 \centering
  \includegraphics[width=0.6\textwidth]{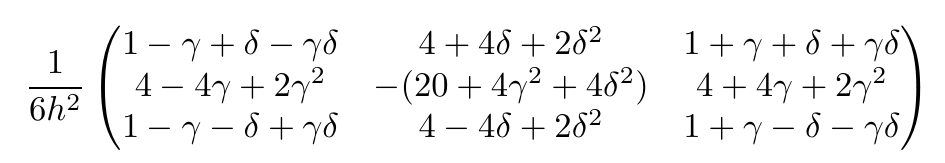}
\caption{The compact nine-point coefficients of the FD scheme from   \cite{MR1440337}.} \label{stencil}
\end{figure}


\section{High order compact schemes for anisotropic elliptic PDEs} \label{sec:anis}

Using the same idea and approach, we have also developed a HOC scheme for anisotropic elliptic partial differential equations
\eqm
   A_{11} u_{xx}   + 2A_{12} u_{xy}  + A_{22} u_{yy} + a u_x + b u_y + Ku = f, \label{anis-PDE}
\enm
with constant coefficients,  and  Dirichlet, Neumann, and Robin boundary conditions. The method is exact for any fourth order polynomials if $K=0$. For well-posedness of the PDE, we assume that $A_{12}^2 - A_{11} A_{22} <0$.
We have found that the HOC scheme works better when  $A_{11} =A_{22}$. If this condition is not true,  we can use a scaling in one coordinate direction to transform the PDE so that $A_{11} =A_{22}$.

A Robin boundary condition has the following form
\eqm
  {\bf A } \grad u \cdot {\bf n} + \sigma u = g.
\enm
For example, if the Robin boundary is defined at $x=x_l$, then ${\bf n} = (-1,0)$ and the BC becomes
\eqm
 \left.  \left (- A_{11} u_x - A_{12}  u_y   + \sigma u \right )  \frac{\null}{\null} \! \right |_{(x_l,y)} = g(y) .
\enm

As  before,  the HOC scheme uses the 9-point stencil at interior grid points  can be written as
\eqmno
&&  \dsp  \sum_{i_k=-1}^{1} \sum_{j_k=-1}^{1} \alpha_{i_k,j_k} U_{i+i_k,j+j_k}  =  \sum_{i_k=-1}^{1} \sum_{j_k=-1}^{1}  \beta_{i_k,j_k} f(x_{i+i_k},y_{j+j_k})  , \\ 
  && \dsp  \sum_{i_k,j_k=-1}^1 \!\!  \beta_{i_k,j_k} = 1.
\enmno
The differences from previous HOC schemes is that the system of linear equations for the coefficients, and the PDE relations are different. We list the linear system of equations and those PDE relations in  Appendix.

 \subsection{FD coefficients for diffusion-advection equations at interior grid points}

There are infinite number of solution sets of the FD scheme. Using the Matlab/Maple symbolic computation package, we have obtained a particular set of the coefficients for {\em interior grid points} below assuming that $A_{11}=A_{22}$ and  zero convection coefficients ($a=b=0$):
\eqml{anis_coeffA}
&U_{i,j}: \dsp  \;   \frac{1}{\lambda h^2}\left[
               \begin{array}{ccc}
                 2(A_{11} - 2 A_{12})(A_{11} - A_{12}) & 8(A_{11}^2-A_{12}^2) & 2(A_{11} + A_{12})(A_{11} + 2A_{12}) \\
                 8(A_{11}^2-A_{12}^2) & \lambda K h^2 -8(5 A_{11}^2 - 2 A_{12}^2) & 8(A_{11}^2-A_{12}^2) \\
                 2(A_{11} + A_{12})(A_{11} + 2A_{12})& 8(A_{11}^2-A_{12}^2) & 2(A_{11} - 2A_{12})(A_{11} - A_{12}) \\
               \end{array}
             \right],
             \enml

             \eqml{anis_coeffB}
&f_{i,j} \dsp:  \;   \frac{1}{2\lambda}\left[
               \begin{array}{ccc}
                 A_{11} - A_{12} & 0 & A_{11} + A_{12} \\
                 0 & 20 A_{11} - 2Kh^2 & 0 \\
                 A_{11} + A_{12} & 0 & A_{11} - A_{12} \\
               \end{array}
             \right],
\enml
where $\lambda = 12A_{11}-Kh^2$. We can see that when $A_{11} > 2 |A_{12}|$  and $K\le 0$, the coefficient matrix of the FD  scheme is an M-matrix.
For non-zero convection coefficients, the expressions for the set of coefficients are complicated, it is easier to use and store numerical solutions.


 \begin{theorem}
 Let $u(x,y)\in C^6({\cal{R}})$ be the solution to \eqn{anis-PDE} with a Dirichlet  boundary condition, $U_{ij}$ be the finite difference solution obtained from the HOC scheme. Then,
 \begin{enumerate}[label=\arabic*), itemindent=0.5em]
 \item the finite difference scheme with a set of coefficients given by (\ref{anis_coeffA})-\eqn{anis_coeffB} is fourth order accurate if $A_{11} > 2 |A_{12}|$, $K\le 0$,  $a=0$, $b=0$;
 \item  the finite difference scheme is asymptotically fourth order convergent if $A_{11} > 2 |A_{12}|$ and $K\le 0$  with  general constants $a$ and $b$.
\end{enumerate}
%
\end{theorem}

{\bf Proof:} For the first case, from the construction of the finite difference coefficients, we know that the local truncation errors are at least $O(h^3)$ since we have matched all terms of fourth order partial derivatives.   We can see that the FD  scheme for $U_{ij}$ has the central symmetry, which means that all the coefficients in odd partial derivatives canceled out in the Taylor expansion of the local truncation errors.  When $A_{11} > 2 |A_{12}|$ and $K\le 0$, the coefficient matrix of the FD equations is an M-matrix. Thus, from the  convergence theorem in \cite{morton-mayers}, we conclude the fourth order convergence.

For  general anisotropic diffusion and advection equations with a Dirichlet boundary condition, from the continuity of the  FD coefficients  and the fact that the  coefficients involving advection terms are of $O(\|{\bf a}\| \|\bfalf\| h^5)\sim O(h^4)$ in the local truncation errors, we can conclude asymptotic fourth order convergence. \hfill $\square$

\subsection{Finite difference coefficients for flux boundary conditions}

If a Robin boundary condition is specified, we need to solve another set of coefficients to take into account the boundary condition. The procedure is the same as before except that the linear system of equations and the  high order PDE relations are different. Unfortunately, the linear system of equations using the same approach is not consistent. Nevertheless, the inconsistency comes from fourth order partial derivatives. Thus, we can simply use the singular value solution (SVD) or ignore some terms  such as $x^4$ and $y^4$  to get a consistent system. With the second approach, we have obtained a set of coefficients below when  $K=0$, $a=0$, $b=0$
 assuming a Robin BC  at the left boundary using the Maple package:
\begin{equation}\label{Coeff_ansi_A}
U_{i,j}:  \; \frac{1}{6 A_{11} h^2} \left[\begin{array}{cc}  - A_{11} \sigma h + (3 A_{11}^2 + 6 A_{11} A_{12} - 4 A_{12}^2)& 3 A_{11}^2 + 6 A_{11} A_{12} + 4 A_{12}^2 \\ - 10 A_{11} \sigma h + (8 A_{12}^2 - 18 A_{11}^2)& 6 A_{11}^2 - 8 A_{12}^2\\- A_{11} \sigma h + (3 A_{11}^2 - 6 A_{11} A_{12} - 4 A_{12}^2) & 3 A_{11}^2 - 6 A_{11} A_{12} + 4 A_{12}^2 \end{array}\right],
\end{equation}
\begin{equation} f_{i,j}:  \;   \frac{1}{48 A_{11}(A_{11}^2 - 2 A_{12}^2)}
\left[\begin{array}{ccc} \lambda+\mu &0 & \zeta+\eta\\
0& 4A_{11} (9 A_{11}^2 - {20 A_{12}^2}) &0 \\
\lambda-\mu& 0& \zeta-\eta \end{array}\right],
\end{equation}
\begin{equation}\label{Coeff_ansi_C}
g_{j}:  \;  \frac{1}{6 h}
\left[\begin{array}{c} - 1 \\  -10 \\ -1 \end{array}\right],
\end{equation}
where

\vspace{-1.cm}
\begin{align*}
  \lambda &= - A_{11}^3  + 4 A_{11} A_{12}^2 ,  \quad \mu = 3 A_{11}^2 A_{12} - 4 A_{12}^3, \\
  \zeta &= 7 A_{11}^3 - 12 A_{11} A_{12}^2, \quad  \eta= 5 A_{11}^2 A_{12} -12 A_{12}^3.
\end{align*}
\noindent
Note that  the coefficient matrix from this set of FD coefficients and \eqn{anis_coeffA}  is an M-matrix if $h$ is
small enough, $K\le 0$, and $A_{11} >  1+\frac{\sqrt{21}}{3} |A_{12}|\approx 2.5275 |A_{12}| $, a stronger condition than that of  interior grid points. We have a  super-third convergence theorem of the finite difference scheme.

\begin{theorem}
 Let $u(x,y)\in C^5({\cal{R}})$ be the solution to \eqn{anis-PDE} with a Robin  boundary condition, $U_{ij}$ be the finite difference solution obtained from the HOC scheme. Then,
 \begin{enumerate}[label=\arabic*), itemindent=0.5em]
 \item the finite difference scheme with a set of coefficients given by (\ref{anis_coeffA})-\eqn{anis_coeffB} at interior grid points, (\ref{Coeff_ansi_A})-\eqn{Coeff_ansi_C} at  boundary grid points,  is super-third  accurate if $A_{11} >  1+\frac{\sqrt{21}}{3}  |A_{12}|$, $K=0$, $a=0$, $b=0$;
 \item  the finite difference scheme is at least asymptotically super-third convergent if $A_{11} >  1+\frac{\sqrt{21}}{3}  |A_{12}|$, $K\le 0$, and $\sigma\ge 0$  with  general constants $a$ and $b$.
\end{enumerate}
\end{theorem}

\ignore{
\begin{theorem}
 Let $u(x,y)\in C^5({\cal{R}})$ be the solution to an anisotropic diffusion equation, $U_{ij}$ be the finite difference solution obtained from the compact scheme with a set of coefficients given by (\ref{Coeff_ansi_A})-\eqn{Coeff_ansi_C}. If $A_{11} \geq  1+\frac{\sqrt{21}}{3} |A_{12}| \approx 2.5275 |A_{12}| $ and  $h$ is
small enough, then we have the following estimate
\eqm
  \left \| u(x_i,y_j)- U_{ij} \right \|_{\infty}  \le C h^{3+},
\enm
for Neumann/Robin  boundary conditions.
\end{theorem}  }


 We  show  two  numerical examples in Table~\ref{tab:anis-diffu-D} using Example~\ref{ex2}, the tougher example with oscillatory solutions.
The second to fifth columns show grid refinement results of the  global  and  local truncation errors, and the convergence orders of Example~\ref{ex2} when the parameters are $A_{11}  = A_{22}  = 3$, $A_{12}  =0.5$,  $K=20$, $k_1 = 3, \,  k_2=15$,  $a=50$, $b=-1$,  and a Dirichlet boundary condition.  We  observe  asymptotic fourth order convergence in the local truncation  as well as the global errors.

\begin{table}[htpb]
\caption{Grid refinement analysis  of the HOC scheme for an {\em anisotropic}  diffusion and advection equation using Example~\ref{ex2}. The results in  columns 2-5 are the global and local truncation errors  and convergence order
of  a Dirichlet BC,  while  columns 6-7 are the global error and order  of a Robin boundary condition.
}
\label{tab:anis-diffu-D}

\vthin

 \begin{center}
\begin{tabular}{|c| c c | | c c|}
  \hline
       $N$   & $\|E\|_{\infty} $  & $order$    &  $\|T_h\|_{\infty}$ & $order$  \\       \hline	
       16	&      2.3762e-03 &  &  1.7131 & \\
32 &	1.6165e-04 &  3.8777 & 	 	1.2748e-01 &   3.7483 \\
64	&  1.0492e-05	&  3.9455  & 8.2910e-03 & 3.9426  \\
128	& 6.6211e-07 & 3.9861 &	5.2320e-04 & 3.9861 \\
256	&   	4.2280e-08 & 3.9690	& 3.3259e-05 &   3.9755 \\
512	& 	2.6740e-09 & 	 3.9975 & 1.9678e-06 &  4.0791 \\
      \hline
\end{tabular}
%
  \begin{tabular}{| c c | }
  \hline
       $\|E\|_{\infty} $  & $order$      \\       \hline	
           1.6756e-05   &     \\
	   1.1523e-06 & 3.8621	 \\
  8.4485e-08 & 3.7697  \\
  	6.8592e-09 & 3.6226  \\
   6.3350e-10	 & 3.4366  \\
 	4.8299e-11 & 3.7133  \\
      \hline
\end{tabular}
\end{center}

\end{table}

In the last two columns of Table~\ref{tab:anis-diffu-D}, we show the result  when a Robin boundary condition is prescribed  at $x=0$. The parameters are $A_{11}  = A_{22}  = 1$, $A_{12}  =0.25$,  $K=-2$, $k_1 = 2,  \,  k_2=1$, and the convection coefficients $a=1$, $b=-5$, and $\sigma=-3$. The results show a super-third convergence order (better than 3.5) but not complete fourth order. It is still ongoing project to see whether it is possible to get a full fourth order convergence  for anisotropic diffusion and advection equations with flux boundary conditions.

\section{The fourth order compact scheme for flux BCs in 3D} \label{sec:3D}

The same idea and methodology have also been applied to  three dimensional (3D) problems with flux type of boundary conditions
\eqm
  &    u_{xx}    +  u_{yy}  + u_{zz}   + K u = f(x,y,z), \qquad \mathbf{x}  \in \cal{R},\\ 
  & \left . \frac{\null}{\null} u \right |_{\partial  { \cal{R}}_1} = u_1(\mathbf{x}), \qquad   \left .  \lp \dsp \frac{\partial u}{\partial n}  +  \sigma  u (\mathbf{x})  \rp \right |_{\partial  {\cal{R}}_2} = g(\mathbf{x}).
\enm
 The main challenge probably is in the implementation and indexing.

Without loss of generality, we assume that domain is a cubic $[x_l,\,x_r]\times [y_l,\,y_r] \times [z_l,\,z_r]$, and a  Robin boundary condition is specified at $x=x_l$, and Dirichlet boundary conditions on other parts of the boundary.  With a uniform mesh,  the 19-point fourth-order compact  FD scheme for a Poisson equation on an interior point $(x_i,y_j,z_k)$ with mesh-size $h$ is given by~\cite{MR1377177}
\eqml{compact}
&&\dsp -\frac{4}{h^2}U_{ijk}+\frac{1}{3h^2}\left(U_{i+1,j,k}+U_{i-1,j,k}+U_{i,j+1,k}+U_{i,j-1,k}+U_{i,j,k+1}+U_{i,j,k-1}\right) \\ 
&&\dsp \null \qquad +\frac{1}{6h^2}(U_{i+1,j+1,k}+U_{i+1,j-1,k}+U_{i-1,j+1,k}+U_{i-1,j-1,k}+U_{i+1,j,k+1}+U_{i+1,j,k-1}\\ \eqsp
&&\dsp \null \qquad  +U_{i-1,j,k+1}+U_{i-1,j,k-1}+U_{i,j+1,k+1}+U_{i,j-1,k+1}+U_{i,j+1,k-1}+U_{i,j-1,k-1}) \\ 
&& \null \qquad  =\dsp \frac{1}{12}(6f_{ijk}+f_{i+1,j,k}+f_{i-1,j,k}+f_{i,j+1,k}+f_{i,j-1,k}+f_{i,j,k+1}+f_{i,j,k-1}).
\enml
The $Ku$ term is treated as a source term. The stencil notation for the approximation is
\ignore{
\begin{align}
 U_{i_c,j_c,k}:  \;  \frac{1}{6h^2}\left[
               \begin{array}{ccc}
                 1 & 2 & 1 \\
                 2 & -24 & 2 \\
                 1 & 2 & 1 \\
               \end{array}
             \right],   \quad U_{i_c, j_c, k\pm 1}: \; \frac{1}{6h^2}\left[
               \begin{array}{ccc}
                 0 & 1 & 0 \\
                 1 & 2 & 1 \\
                 0 & 1 & 0 \\
               \end{array}
             \right ] , mbox{TEST}
\end{align}
\begin{align}
 f_{i_c,j_c,k}:  \;  \frac{1}{12}\left[
               \begin{array}{ccc}
                 0 & 1 & 0 \\
                 1 & 6 & 1 \\
                 0 & 1 & 0 \\
               \end{array}
             \right],   \qquad f_{i_c, j_c, k\pm 1}: \;   \frac{1}{12}\left[
               \begin{array}{ccc}
                 0 & 0 & 0 \\
                 0 & 1 & 0 \\
                 0 & 0 & 0 \\
               \end{array}
             \right],
\end{align}
where $i_c=i -1, \, i, \, i +1$ and $j_c=j -1, \, j, \, j +1$. Note that there are seven $f_{ijk}$'s involved in the scheme.}
%
%
  \begin{align}
 U_{i_c,j,k}:  \quad  \frac{1}{6h^2}\left[
               \begin{array}{ccc}
                 1 & 2 & 1 \\
                 2 & -24 & 2 \\
                 1 & 2 & 1 \\
               \end{array}
             \right], \quad  \qquad U_{i_c \pm 1, j, k}: \quad \frac{1}{6h^2}\left[
               \begin{array}{ccc}
                 0 & 1 & 0 \\
                 1 & 2 & 1 \\
                 0 & 1 & 0 \\
               \end{array}
             \right],
\end{align}
%
\begin{align}
 f_{i_c,j,k}:  \quad  \frac{1}{12}\left[
               \begin{array}{ccc}
                 0 & 1 & 0 \\
                 1 & 6 & 1 \\
                 0 & 1 & 0 \\
               \end{array}
             \right], \quad  \qquad f_{i_c \pm 1, j, k}: \quad   \frac{1}{12}\left[
               \begin{array}{ccc}
                 0 & 0 & 0 \\
                 0 & 1 & 0 \\
                 0 & 0 & 0 \\
               \end{array}
             \right],
\end{align}
where for example, $j=j_c -1, \, j_c, \, j_c +1$ and $k=k_c -1, \, k_c, \, k_c +1$.
Note that there are seven $f_{ijk}$'s in the scheme.

For a grid point $(x_0, y_j, z_k)$ at the boundary $x=x_q$ where a   Robin boundary condition is prescribed, similar to the 2D case, the fourth order FD equation has the following  form,
\eqm
&&   \dsp  \sum_{i_l=0}^{1} \sum_{j_l=-1}^{1}  \sum_{k_l=-1}^{1} \!\!\ \alpha_{i_l,j_l,k_l} U_{i_l,j+j_l,k+k_l}   =  \sum_{i_l=-1}^{1} \sum_{j_l=-1}^{1}  \sum_{k_l=-1}^{1}  \!\!\  \beta_{i_l,j_l,k_l}  f(x_{i_l},y_{j+j_l},z_{k+k_l})  \nonumber \\ 
&& \dsp  \null \hspace{3cm} \qquad \null \qquad \null +\,   \sum_{j_l=-1}^{1}  \sum_{k_l=-1}^{1} \!\!\
 \gamma_{j_l,k_l}    \, g(y_{j+j_l},z_{k+k_l})  ,  \\  
 && \dsp  \null \hspace{1cm} \qquad   \dsp  \sum_{i_l=-1}^{1} \sum_{j_l=-1}^{1} \sum_{k_l=-1}^{1}  \!\!\ \beta_{i_l,j_l,k_l}  = 1.
\enm
Thus, if all points are involved, then the total degree of freedom is $18+27+9=54$. We want the FD method is exact (when $K=\sigma=0$) if the solution is any fourth order polynomials $\dsp \sum_{0\le i+j+k\le 4} \!\!\!\!\!\!x^i y^j z^k$ that would lead to $35$ equations in addition to the constraint. Thus, we  have an under-determined linear system of equations. To enforce the discrete maximum principle, we enforce the sign property so that the coefficient matrix for the finite difference scheme is an M-matrix.
In general, we have infinite number of solutions and we can pick up good ones that can have expected symmetries and the least non-zero coefficients. For a Poisson equation with a  Robin boundary condition at $x=x_l$, once again, we have a  simple analytic solution set as listed below:
 \begin{align}
U_{0,j,k}:  \;   \frac{1}{6h^2}\left[
               \begin{array}{ccc}
                 1 & 2 & 1 \\
                 2 & -12(2+\sigma) & 2 \\
                 1 & 2 & 1 \\
               \end{array}
             \right], \qquad   U_{1,j,k}:  \;  \frac{1}{6h^2}\left[
               \begin{array}{ccc}
                 0 & 2 & 0 \\
                 2 & 4 & 2 \\
                 0 & 2 & 0 \\
               \end{array}
             \right],\label{3Dcompact1}
\end{align}
\begin{align}
 f_{0 ,j,k}: \; \frac{1}{12}\left[
               \begin{array}{ccc}
                 0 & 1 & 0 \\
                 1 & 6 & 1 \\
                 0 & 1 & 0 \\
               \end{array}
             \right],   \quad  f_{-1,j,k}: \;  \frac{1}{12}\left[
               \begin{array}{ccc}
                 0 & 0 & 0 \\
                 0 & -1 & 0 \\
                 0 & 0 & 0 \\
               \end{array}
             \right],
             \quad  f_{1,j,k}: \;  \frac{1}{12}\left[
               \begin{array}{ccc}
                 0 & 0 & 0 \\
                 0 & 3 & 0 \\
                 0 & 0 & 0 \\
               \end{array}
             \right],
\end{align}
\begin{align}
 g_{j,k}: \;   \frac{1}{h}\left[
               \begin{array}{ccc}
                 0 & 0 & 0 \\
                 0 & -2 & 0 \\
                 0 & 0 & 0 \\
               \end{array}
             \right].\label{3Dcompact3}
\end{align}
Note that the discretization of the Neumann boundary condition is the same as the ghost point method, but the right hand side $f_{0jk}$ has been adjusted.

\begin{theorem} \label{thm-3D}
Let $U_{ijk}$ be the finite difference solution obtained from the fourth order scheme with the set of the solution of the finite difference coefficients given above. 
Then, the algorithm is {\em exact}  if the solution is any fourth  order polynomials  when $K=\sigma=0$. For general solutions
  $u(x,y)\in C^6(\cal{R})$,   we have the following error estimate,
\eqm
  \left \| u(x_i,y_j,z_k)- U_{ijk} \right \|_{\infty}  \le C h^4.
\enm
\end{theorem}

\begin{example} \label{ex3}
A three dimensional example:
\eqmno 
 &&  \dsp u(x,y,z) = \sin(k_1 x) \sin(k_2 y) \sin(k_3 z), \quad \quad  (x,y,z) \in (0,\;1)^3, \\ 
 &&  \dsp f(x,y,z) =  - \left (k_1^2 + k_2^2 + k_3^2 \right )   \sin(k_1 x)  \sin(k_2 y) \sin(k_3 z) ,  \quad (x,y,z) \in (0,\;1)^3, \\ 
 && \dsp  \left. \lp  \frac{\partial u}{\partial n} + \sigma u \rp  \right |_{x=0} =  \left  (-k_1\; \cos(k_1 x) + \sigma  \sin(k_1 x)
\frac{\null}{\null}  \right )  \sin(k_2 y) \sin(k_3 z) ,\qquad (y,z) \in  (0,\;1)^2  .
 \enmno
\end{example}

In Table~\ref{tab:3D}, we show a grid refinement analysis with $k_1=1, k_2=2, k_3=10$ and $\sigma=2$. Note that, the total degree of freedom of the finite difference equations when $N=128$ is $2,064,512$, more than two millions.
We adopt an extrapolation cascadic multigrid method with the conjugate gradient (CG) smoother (or BiCGStab if $A_h$ is non-symmetric) from \cite{pan2017extrapolation} to solve the sparse linear system.
The CPU time is recorded  on  a Laptop with an Intel(R) Core(TM) i7-8565U CPU @ 1.80GHz and 8.0 GB RAM.
We can  see clearly a fourth order convergence.

\begin{table}[htpb]
\caption{A grid refinement analysis for the 3D example with a Robin BC at $x=0$. Fourth order convergence can be seen clearly.
\textcolor{blue}{An extrapolation cascadic multigrid method with the conjugate gradient (CG) smoother \cite{pan2017extrapolation} is used to solve the linear systems. }
} \label{tab:3D}
 \begin{center}
\begin{tabular}{||c| c c   || r ||}
  \hline
       $N$   & $\|E\|_{\infty}$  & $order$  &       $CPU$       \\       \hline	
    8 & 5.9285e-03  &                        &       $<$0.01 s     \\
   16 & 3.8247e-04  & 3.9543                 &        0.02 s     \\
   32 & 2.3691e-05  & 4.0129                 &               0.13 s     \\
   64 & 1.4883e-06  & 3.9926                 &          1.32 s     \\
  128 & 9.2985e-08  & 4.0006                 &         11.00 s    \\
\hline
\end{tabular} 
 \end{center}
\end{table}

\section{Conclusions and discussions}  \label{sec:conclud}

In this paper, we have solved an important problem in computational mathematics, that is, whether there  exist fourth order compact schemes for Poisson, Helmholtz, and diffusion-advection equations with flux type boundary conditions.
The answer is yes if we can extend the source term $f$ to one grid line (surface in 3D) with a quadratic extension that is third  order accurate. Without the $f$-extension, then we probably can only achieve super-third  convergence in which the HOC methods have  been developed in this paper.
Using a brand new approach, we have developed  new  fourth order compact schemes in both 2D and 3D that can guarantee the consistence, stability, so the convergence.

The new idea and methodology have also been applied  to anisotropic diffusion and advection equations with Dirichlet, Neumann, or Robin BCs with constant coefficients. Fourth order convergence has been proved for Dirichlet boundary conditions while super-third convergence has been proved for flux boundary conditions. So it is still an open question to develop fourth order  compact schemes for anisotropic diffusion and advection equations with flux boundary conditions.
Technically, the new idea can be applied to BVP of elliptic PDEs with variable coefficients but are not recommended due to the  computational  cost.  This is because  for variable coefficient PDEs, the coefficient matrix  for the finite difference and weight coefficients  would be different at every grid point. We need to solve $O(N^2)$ such systems of equations with the sign constraint for 2D problems. Thus, the computational cost would be overwhelming. We also need the first, second order partial derivatives of the variable coefficients.
For constant coefficient PDE, the coefficients of the finite difference equations are just needed to be computed once for all. 

While the developed methods are for rectangular domains, 
we think the methods still can work for polygonal domain such as an L-shaped domain as long as all boundary segments are parallel to one of the axes. For general curved boundary, we can derive third order accurate compact schemes using the augmented approach developed in \cite{pan-he-li-HOC21}. But we are not sure about fourth order compact schemes and think it will be very challenging for developing such schemes.

%
%



\parskip= 0.0cm
\renewcommand{\baselinestretch}{1.0}

\bibliographystyle{amsplain}

 \bibliography{hoc,../../TEX/BIB/bib,../../TEX/BIB/zhilin,../../TEX/BIB/other}

\appendix


\vthin

\section{High order PDE relations for general elliptic PDEs with constant coefficients}
To derive the fourth order compact finite difference scheme for an anisotropic  diffusion and advection equations,
\eqml{anis-PDEB}
  \dsp    A_{11} u_{xx}   + 2A_{12} u_{xy}  + A_{22} u_{yy} + a u_x + b u_y + Ku = f,
\enml
with a Dirichlet or Robin boundary conditions, we need more PDE relations in order to derive high order compact schemes.

\begin{lemma}
Let $u(x,y)\in C^5({\cal R})$ be the  solution to \eqn{anis-PDEB}, then the following relations are true.
\end{lemma}

\vspace{-0.6cm}
\eqmno
  && A_{11} u_{xxx}   + 2A_{12} u_{xxy}  + A_{22} u_{xyy} + a u_{xx} + b u_{xy} + Ku_x = f_x\\ 
  && A_{11} u_{xxy}   + 2A_{12} u_{xyy}  + A_{22} u_{yyy} + a u_{xy} + b u_{yy} + Ku_y = f_y\\ 
  && A_{11} u_{xxxx}   + 2A_{12} u_{xxxy}  + A_{22} u_{yyxx} + a u_{xxx} + b u_{xxy} + Ku_{xx} = f_{xx}\\ 
  && A_{11} u_{xxxy}   + 2A_{12} u_{xxyy}  + A_{22} u_{xyyy} + a u_{xxy} + b u_{xyy} + Ku_{xy} = f_{xy}\\ 
  && A_{11} u_{xxyy}   + 2A_{12} u_{xyyy}  + A_{22} u_{yyyy} + a u_{xyy} + b u_{yyy} + Ku_{yy} = f_{yy}.
\enmno

\subsection{The linear system of equations of the FD coefficients at interior grid points for diffusion and advection equations}
The finite difference coefficients
 $\alpha_{i_k,j_k}$ and $\beta_{i_k,j_k}$ for interior grid points  are determined from the following system of equations.  The first six equations (required for all quadratic polynomials) are,
\eqml{anisA}
&&\dsp  \sum_{i_k=-1}^{1} \sum_{j_k=-1}^{1} \alpha_{i_k,j_k}    -  K = 0\\  
&&\dsp  \sum_{i_k=-1}^{1} \sum_{j_k=-1}^{1} \alpha_{i_k,j_k}  h_{i_k}  -  K \!\!  \sum_{i_k=-1}^{1} \sum_{j_k=-1}^{1}  \beta_{i_k,j_k} h_{i_k}   =0  \\ 
&&\dsp  \sum_{i_k=-1}^{1} \sum_{j_k=-1}^{1} \alpha_{i_k,j_k}  h_{j_k}  -  K \!\!  \sum_{i_k=-1}^{1} \sum_{j_k=-1}^{1}  \beta_{i_k,j_k} h_{j_k}   =0  \\ 
&& \dsp \sum_{i_k=-1}^{1} \sum_{j_k=-1}^{1} \alpha_{i_k,j_k}  \frac{h_{i_k}^2 } {2} -  \sum_{i_k=-1}^{1} \sum_{j_k=-1}^{1}  \beta_{i_k,j_k}   \left ( A_{11} + \frac{ h_{i_k}^2 } {2}  K + a h_{i_k} \right ) = 0 \\ 
&& \dsp \sum_{i_k=-1}^{1} \sum_{j_k=-1}^{1} \alpha_{i_k,j_k}  \frac{h_{j_k}^2 } {2} -  \sum_{i_k=-1}^{1} \sum_{j_k=-1}^{1}   \beta_{i_k,j_k}  \left (A_{22}+ \frac{ h_{j_k}^2 } {2}  K  + b  h_{j_k} \right ) = 0 \\ 
 && \dsp \sum_{i_k=-1}^{1} \sum_{j_k=-1}^{1} \alpha_{i_k,j_k} h_{i_k} h_{j_k} -   \sum_{i_k=-1}^{1} \sum_{j_k=-1}^{1}  \beta_{i_k,j_k}  \left (  2 A_{12} +
 h_{i_k} h_{j_k} K + b  h_{i_k} + a  h_{j_k} \right ) = 0.
\enml
The next four equations (required for cubic polynomials) are
\eqml{anis-coeB-diffu}
&& \dsp \sum_{i_k=-1}^{1} \sum_{j_k=-1}^{1} \alpha_{i_k,j_k}  \frac{h_{i_k}^3 } {3!} -  \sum_{i_k=-1}^{1} \sum_{j_k=-1}^{1}  \beta_{i_k,j_k} h_{i_k} \left ( A_{11}+  a\, \frac{ h_{i_k} } {2}   \right ) = 0 \\ 
&& \dsp \sum_{i_k=0}^{1} \sum_{j_k=-1}^{1} \alpha_{i_k,j_k}  \frac{h_{i_k}^2  h_{j_k}} {2} -  \sum_{i_k=-1}^{1} \sum_{j_k=-1}^{1}  \beta_{i_k,j_k} \left ( A_{11} h_{j_k}  + 2A_{12 }   h_{i_k} + a \,  h_{i_k} h_{j_k} + b \, \frac{ h_{i_k}^2 } {2} \right ) = 0 \\ 
&& \dsp \sum_{i_k=-1}^{1} \sum_{j_k=-1}^{1} \alpha_{i_k,j_k}  \frac{h_{i_k}  h_{j_k}^2} {2} -  \sum_{i_k=-1}^{1} \sum_{j_k=-1}^{1}  \beta_{i_k,j_k} \left ( A_{22} h_{i_k}  +  2A_{12} h_{j_k}   + a \,   \frac{ h_{j_k}^2 } {2} + b\, h_{i_k} h_{j_k}  \right ) = 0 \\ 
&&\dsp  \sum_{i_k=-1}^{1} \sum_{j_k=-1}^{1} \alpha_{i_k,j_k}  \frac{  h_{j_k}^3} {3!} -  \sum_{i_k=-1}^{1} \sum_{j_k=-1}^{1}  \beta_{i_k,j_k} h_{j_k} \left ( A_{22} +  b\, \frac{ h_{j_k} } {2}   \right ) = 0  .
\enml

The next  five equations (required for quartic polynomials) are
 \eqml{anis-coeC-diffu}
&& \dsp \sum_{i_k=-1}^{1} \sum_{j_k=-1}^{1} \alpha_{i_k,j_k}  \frac{h_{i_k}^4 } {4!} -  \sum_{i_k=-1}^{1} \sum_{j_k=-1}^{1}  \beta_{i_k,j_k}  \frac{h_{i_k}^2  } {2} A_{11}  = 0 \\ \eqsp
&& \dsp \sum_{i_k=-1}^{1} \sum_{j_k=-1}^{1} \alpha_{i_k,j_k}  \frac{h_{i_k}^3  h_{j_k}} {3!} -  \sum_{i_k=-1}^{1} \sum_{j_k=-1}^{1}  \beta_{i_k,j_k} \left ( h_{i_k} h_{j_k}  A_{11} +  {h_{i_k}^2 }  A_{12}   \right ) = 0  \\ \eqsp
&& \dsp \sum_{i_k=-1}^{1} \sum_{j_k=-1}^{1} \alpha_{i_k,j_k}  \frac{6 h_{i_k}^2 h_{j_k}^2} {4!} -  \sum_{i_k=-1}^{1} \sum_{j_k=-1}^{1}  \beta_{i_k,j_k} \left (   \frac{h_{j_k}^2  } {2} A_{11} +  2h_{i_k} h_{j_k}A_{12}  +  \frac{h_{i_k}^2  } {2} A_{22} \right )  = 0 \\ \eqsp
&& \dsp \sum_{i_k=-1}^{1} \sum_{j_k=-1}^{1} \alpha_{i_k,j_k}  \frac{h_{i_k}  h_{j_k}^3} {3!} -  \sum_{i_k=-1}^{1} \sum_{j_k=-1}^{1}  \beta_{i_k,j_k}  \left ( h_{i_k} h_{j_k}  A_{22} +  {h_{j_k}^2 } A_{12}   \right )   = 0 \\ \eqsp
&& \dsp \sum_{i_k=-1}^{1} \sum_{j_k=-1}^{1} \alpha_{i_k,j_k}  \frac{  h_{j_k}^4} {4!} -  \sum_{i_k=-1}^{1} \sum_{j_k=-1}^{1}  \beta_{i_k,j_k} h_{j_k} \frac{h_{j_k}^2 } {2} A_{22} = 0 . 
 \enml

\subsection{The linear system of equations of the FD coefficients at a flux boundary grid points for diffusion and advection equations}

For a flux Robin condition. The linear  system of equations for the FD coefficients  has similar form as above but the summation $\dsp \sum_{i_k=-1}^{1} \left( \null \frac{\null}{\null} \right )$ changes to $\dsp \sum_{i_k=0}^{1}\left( \null \frac{\null}{\null} \right )$ in all equations involving $\alpha_{i_k,j_k}$ terms.
The first-third, fifth-sixth,  nine-tenth, fourteenth,  and sixteenth equations will need to be changed to
\eqml{anis-BC-eq}
&&\dsp  \sum_{i_k=0}^{1} \sum_{j_k=-1}^{1} \alpha_{i_k,j_k}    -  K -   \sigma  \sum_{j_k=-1}^1 \gamma_{j_k} = 0\\ \eqsp
&&\dsp  \sum_{i_k=0}^{1} \sum_{j_k=-1}^{1} \alpha_{i_k,j_k}  h_{i_k}  -  K \!\!  \sum_{i_k=-1}^{1} \sum_{j_k=-1}^{1}  \beta_{i_k,j_k} h_{i_k}  +  \!\! \sum_{j_k=-1}^1 \!\! \gamma_{j_k}   A_{11} =0  \\ \eqsp
  &&\dsp  \sum_{i_k=0}^{1} \sum_{j_k=-1}^{1} \alpha_{i_k,j_k}  h_{j_k}  -  K \!\!  \sum_{i_k=-1}^{1} \sum_{j_k=-1}^{1}  \beta_{i_k,j_k} h_{j_k}  -    \!\! \sum_{j_k=-1}^1 \!\! \gamma_{j_k}  \lp  \sigma  h_{j_k} -A_{12}\rp  =0  \\ \eqsp
&& \dsp \sum_{i_k=0}^{1} \sum_{j_k=-1}^{1} \alpha_{i_k,j_k}  \frac{h_{j_k}^2 } {2} -  \sum_{i_k=-1}^{1} \sum_{j_k=-1}^{1}   \beta_{i_k,j_k}  \left (A_{22}+ \frac{ h_{j_k}^2 } {2}  K  + b  h_{j_k} \right ) -  \sigma  \!\! \sum_{j_k=-1}^1 \!\! \gamma_{j_k}  \lp \frac{h_{j_k}^2 } {2}   -    A_{12}  h_{j_k}  \rp = 0 \\ \eqsp
 && \dsp \sum_{i_k=0}^{1} \sum_{j_k=-1}^{1} \alpha_{i_k,j_k} h_{i_k} h_{j_k} -   \sum_{i_k=-1}^{1} \sum_{j_k=-1}^{1}  \beta_{i_k,j_k}  \left (  2 A_{12} + h_{i_k} h_{j_k} K + b  h_{i_k} + a  h_{j_k} \right ) +  \!\! \sum_{j_k=-1}^1 \!\! \gamma_{j_k}   A_{11}  h_{j_k}  = 0 \\ \eqsp
 &&\dsp  \sum_{i_k=0}^{1} \sum_{j_k=-1}^{1} \alpha_{i_k,j_k}  \frac{  h_{j_k}^3} {3!} -  \sum_{i_k=-1}^{1} \sum_{j_k=-1}^{1}  \beta_{i_k,j_k} h_{j_k} \left ( A_{22} +  b\, \frac{ h_{j_k} } {2}   \right ) -  \sigma \sum_{j_k=-1}^1 \gamma_{j_k}   \frac{h_{j_k}^3 } {3!}= 0  \\ \eqsp
 && \dsp \sum_{i_k=0}^{1} \sum_{j_k=-1}^{1} \alpha_{i_k,j_k}  \frac{h_{i_k}  h_{j_k}^2} {2} -  \sum_{i_k=-1}^{1} \sum_{j_k=-1}^{1}  \beta_{i_k,j_k} \left ( A_{22} h_{i_k}  +  2A_{12} h_{j_k}   + a \,   \frac{ h_{j_k}^2 } {2} + b\, h_{i_k} h_{j_k}  \right )   \\ \eqsp
  && \dsp \null \qquad \qquad \hspace{5cm} \null  -  \!\! \sum_{j_k=-1}^1 \!\! \gamma_{j_k}  \lp  \sigma  \frac{h_{j_k}^3 } {3!}
 -    A_{12}  \frac{h_{j_k}^2 } {2} \rp= 0 \\ \eqsp
 && \dsp \sum_{i_k=0}^{1} \sum_{j_k=-1}^{1} \alpha_{i_k,j_k}  \frac{h_{i_k}  h_{j_k}^3} {3!} -  \sum_{i_k=-1}^{1} \sum_{j_k=0}^{1}  \beta_{i_k,j_k}  \left ( h_{i_k} h_{j_k}  A_{22} +   {h_{j_k}^2} A_{12}   \right )
 +   \!\! \sum_{j_k=-1}^1 \!\! \gamma_{j_k}   A_{11}   \frac{h_{j_k}^3 } {3!}    = 0 \\ \eqsp
 && \dsp \sum_{i_k=0}^{1} \sum_{j_k=-1}^{1} \alpha_{i_k,j_k}  \frac{  h_{j_k}^4} {4!} -  \sum_{i_k=-1}^{1} \sum_{j_k=-1}^{1}  \beta_{i_k,j_k} h_{j_k} \frac{h_{j_k}^2 } {2} A_{22} +   \!\! \sum_{j_k=-1}^1 \!\! \gamma_{j_k}   A_{12}   \frac{h_{j_k}^3 } {3!}   = 0.
\enml

\subsection{Finite difference coefficients of  $U_{ij}$  for diffusion and advection equations}
\begin{figure}[phtb]
\includegraphics[width=1.1\textwidth]{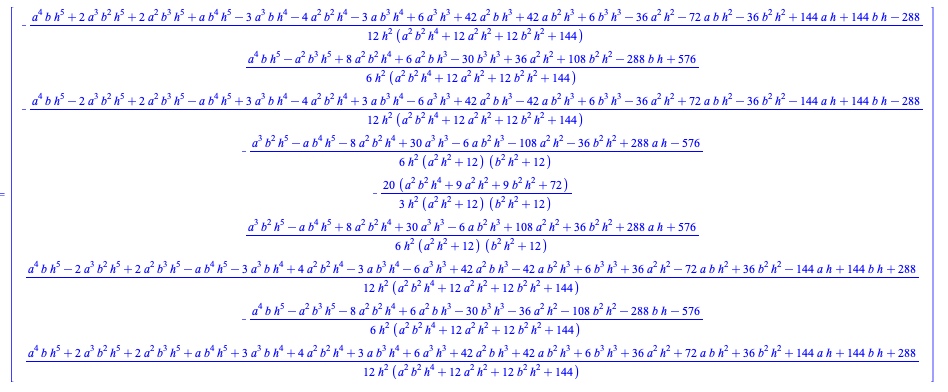}
\caption{Finite difference coefficients for diffusion-advection equations at interior grid points computed with Maple.  Note that the sign property holds when $h$ or $a$ and $b$ are sufficient small.}
\label{FD_Coef_MapleB}
\end{figure}

%
\subsection{Combination coefficients of $f_{ij}$ terms for diffusion and advection equations}
$\null$
\begin{figure}[phbt]
\centerline{\includegraphics[width=0.7\textwidth]{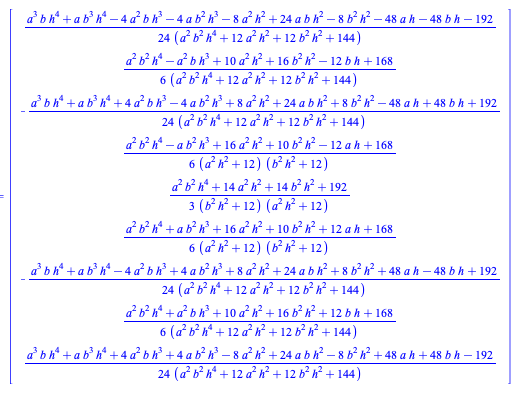}}
\caption{Combination coefficients of $f_{ij}$'s whose sum equals  one for interior   diffusion-advection equations computed at interior grid points with Maple.}
\label{FD_Coef_MapleC}
\end{figure}

\end{document}

%% file: zmacros.tex
\newcommand{\lp}{\left (}
\newcommand{\rp}{\right )}

\newcommand{\eqn}[1]{(\ref{#1})}

\newcommand{\eq}{\begin{equation}}
\newcommand{\en}{\end{equation}}
\newcommand{\eqm}{\begin{eqnarray}}
\newcommand{\enm}{\end{eqnarray}}
\newcommand{\eqmno}{\begin{eqnarray*}}
\newcommand{\enmno}{\end{eqnarray*}}
\newcommand{\eql}[1]{\begin{equation}\label{#1}}
\newcommand{\half}{{\frac{1}{2}}}

\newcommand{\vthin}{\vskip 5pt}

\newcommand{\ignore}[1]{}
\def\bc{\begin{center}}
\def\ec{\end{center}}
\def\bi{\begin{itemize}}
\def\ei{\end{itemize}}
\def\be{\begin{enumerate}}
\def\ee{\end{enumerate}}

\def\alf{\alpha}

\def\reals{{{\rm l} \kern -.15em {\rm R} }}

\def\grad{\nabla}
\def\qquad{\quad\quad}

\def\grad{\nabla}

\def\sf{{\cal F}}

\def\eqsp{\noalign{\vskip 5pt}}
\def\qquad{\quad\quad}

\def\dsp{\displaystyle}
\newcommand{\eqml}[1]{\eql{#1}\begin{array}{rcl}}
\newcommand{\enml}{\end{array}\end{equation}}